\documentstyle[amscd]{amsart}

\def\AA{{\blb A}}

\def\NN{{\blb N}}

\def\PP{{\blb P}}

\def\ZZ{{\blb Z}}

\def\blb#1{\Bbb#1}

\def\11{{1\kern-3.5pt 1}}
\def\mumu{{\mu\kern-4.2pt\mu}}

\def\boxtimes{\setbox0\hbox{$\Box$}\copy0\kern-\wd0\hbox{$\times$}}

\def\Ab{{{\sf Ab}}}

\def\coker{\operatorname {coker}}

\def\hd{\operatorname {hd}}
\def\Hom{\operatorname {Hom}}
\def\im{\operatorname {im}}

\def\ker{\operatorname {ker}}
\def\Ker{\operatorname {ker}}

\def\Spec{\operatorname {Spec}}

\def\th{\operatorname {th}}    

\def\Tor{\operatorname {Tor}}
\def\cTor{\operatorname {{\cal T}\!{\it or}}}

\def\BiGr{{\sf BiGr}}

\def\BiMod{{\sf BiMod}}
\def\Bimod{{\sf Bimod}}
\def\BIMOD{{\sf BIMOD}}

\def\dim{\operatorname{dim}}

\def\Ext{\operatorname{Ext}}

\def\fl.dim{\operatorname{flat.dim}}

\def\grmod{{\sf grmod}}

\def\GrMod{{\sf GrMod}}

\def\Hom{\operatorname{Hom}}

\def\id{\operatorname{id}}
\def\Id{\operatorname{Id}}
\def\Im{\operatorname{Im}}

\def\Inj{{\sf Inj}}

\def\Ker{\operatorname{Ker}}

\def\mod{{\sf mod}}
\def\Mod{{\sf Mod}}

\def\Pic{\operatorname{Pic}}

\def\Proj{{\sf Proj}}

\def\Projnc{{\sf Proj_{nc}}}

\def\rank{\operatorname{rank}}

\def\Spec{\operatorname{Spec}}

\def\sup{\operatorname{sup}}
\def\Supp{\operatorname{Supp}}
\def\SSupp{\operatorname{SSupp}}

\def\tors{{\sf tors}}
\def\Tors{{\sf Tors}}

\def\uTor{\operatorname{\underline{Tor}}}
\def\cHom{\cH om }
\def\cExt{\cE xt }
\def\cTor{\cT or }

\def\ucHom{\underline {\cHom }}

\def\ucTor{\underline {\cTor }}

\def\G{\mathop{\underline{\underline{\it \Gamma}}}\nolimits}

\def\d{\downarrow}

\let\oldtext\text
\def\text#1{\oldtext{\normalshape #1}}

\def\d{\delta}

\def\ve{\varepsilon}

\def\o{\omega}

\def\G{\Gamma}

\def\fm{{\frak m}}

\def\sC{{\sf C}}
\def\sD{{\sf D}}
\def\sE{{\sf E}}

\def\sK{{\sf K}}
\def\sL{{\sf L}}
\def\sM{{\sf M}}

\def\coh{{\sf{ coh}}}

\def\Qcoh{{\sf{ Qcoh}}}

\def\psiol{\bar{\psi}}

\def\cB{{\cal B}}

\def\cD{{\cal D}}
\def\cE{{\cal E}}
\def\cF{{\cal F}}
\def\cG{{\cal G}}

\def\cHom{{\cal H}{\it om}}

\def\cL{{\cal L}}
\def\cM{{\cal M}}
\def\cN{{\cal N}}
\def\cO{{\cal O}}

\def\cQ{{\cal Q}}

\def\cT{{\cal T}}
\def\cV{{\cal V}}
\def\cW{{\cal W}}

\newtheorem{lemma}{Lemma}[section]
\newtheorem{proposition}[lemma]{Proposition}
\newtheorem{theorem}[lemma]{Theorem}
\newtheorem{corollary}[lemma]{Corollary}

{

}

\theoremstyle{definition}

\newtheorem{definition}[lemma]{\sl Definition}

\theoremstyle{remark}

\newtheorem{remark}[lemma]{Remark}

\numberwithin{equation}{section}

\begin{document}

\pagenumbering{arabic}

\title{The Grothendieck group of a quantum projective space bundle}

\author{Izuru Mori and S. Paul Smith}

\address{Department of Mathematics, Shizuoka University, Shizuoka 422-8529, JAPAN}
\email{simouri@@ipc.shizuoka.ac.jp}

\address{Department of Mathematics,  Box 354350,
University of
Washington, Seattle, WA 98195}
\email{smith@@math.washington.edu}

\thanks{S. P. Smith was supported by NSF grants DMS 9701578 and 0070560}

\keywords{}
\subjclass{}

\begin{abstract}
We compute the Grothendieck group of non-commutative
analogues of projective space bundles. Our results specialize to
give the  Grothendieck groups of non-commutative analogues of projective
spaces, and specialize to recover the  Grothendieck group of a usual
projective space bundle over a regular noetherian separated
scheme.  As an application, we develop an intersection theory for
quantum ruled surfaces.
\end{abstract}

\maketitle

\setcounter{section}{-1}
\section{Introduction}

A basic result in the K-theory of algebraic varieties is the
computation of the K-groups for projective space bundles
\cite{Berth}, \cite{Man}, \cite{Q}. In this paper we compute $K_0$
for non-commutative analogues of projective space bundles.

Our motivation comes from the problem of classifying
non-commutative surfaces. This program is in its early stages but,
in analogy with the commutative case,
the non-commutative analogues of ruled surfaces appear to play a
central role \cite{A}, \cite{Pat}, \cite {Vn}.  An
intersection theory will be an essential ingredient in the study of
non-commutative surfaces. One may develop an
intersection theory by defining an intersection multiplicity as a
bilinear $\ZZ$-valued form, the Euler form,
on the Grothendieck group of the surface
(cf. \cite{Jorg} and \cite{MS}).
One of our main results yields a formula (\ref{eq.Krs}) for the
Grothendieck group for quantum ruled surfaces, and we use
this to show that the associated intersection theory gives natural
analogues of the commutative results: if $X$ is a
smooth (commutative!) projective curve and $f:\PP(\cE) \to X$ a quantum
ruled surface over $X$, then fibers do not meet and a section meets
a fiber exactly once (Theorem \ref{thm.int.thy}).

\bigskip \noindent {\it Acknowledgements}
We thank A. Nyman and M. Van den Bergh for showing us their preprints.  We also thank
C. Weibel and an anonymous referee for their constructive criticisms and remarks, and
for pointing out an error in an earlier version of this paper.

\section{Terminology and notation}

The utility of co-monads in non-commutative algebraic geometry was first realized by
Rosenberg \cite{Ros}. Later Van den Bergh
\cite{vdB2} developed very effective methods based on this
idea and replaced the language of
co-monads by terminology and notation that is closer to the
classical language of algebra and algebraic geometry.
We now recall some of  his language \cite{vdB2}.

\subsection{Bimodules \cite[Sect. 3.1]{vdB2}}
\label{sect.bimod}

Let $\sL$ and $\sK$ be abelian categories. The category  $\BIMOD(\sK,\sL)$ of
{\sf weak $\sK$-$\sL$-bimodules} is the opposite of
the category of left exact functors $\sL \to \sK$.
If $\cF={}_\sK\cF_\sL$ is a weak $\sK$-$\sL$-bimodule, we write
$\cHom_\sL(\cF,-)$ for the corresponding left exact functor.  We call
$\cF$ a {\sf bimodule} if $\cHom_\sL(\cF,-)$ has a left adjoint (which we
denote by $- \otimes_\sK \cF$).

The Yoneda embedding $\sK \to \BIMOD(\Ab,\sK)$ realizes
objects in $\sK$ as bimodules and every $\Ab$-$\sK$-bimodule is isomorphic to
an object of $\sK$ in this sense.

The composition of left exact functors is left exact. The composition of  the weak bimoduless (i.e., left exact functors) ${}_\sK\cF_\sL$ and ${}_\sL\cG_\sM$ is denoted by $\cF \otimes_\sL \cG$. This composition is a weak $\sK$-$\sM$-bimodule and
 $\cHom_{\sM}(\cF \otimes_\sL \cG,-)
=\cHom_{\sL}(\cF,\cHom_{\sM}(\cG,-))$.
Thus $\BIMOD(\sK,\sK)$ is a monoidal category with identity
object $o_{\sK}$ defined by $\cHom_{\sK}(o_K,-)=\Id_{\sK}$.

\subsection{Algebras \cite[Sect. 3.1]{vdB2}}
\label{sect.alg}

An
algebra object in  $\BIMOD(\sK,\sK)$ is called a  {\sf weak $\sK$-algebra}:  it is a
triple $(A,\mu,\eta)$, where $A \in \BIMOD(\sK,\sK),$
and $\mu:A \otimes_\sK A \to A$ and $\eta:o_\sK \to A$ are
morphisms of weak bimodules such that the obvious diagrams commute.
If $\cHom_{\sK}(A,-)$ has a left adjoint we drop the adjective
``weak'' and simply call $A$ a $\sK$-algebra.
In this paper we are interested in algebras, not weak algebras.

An {\sf $A$-module} is an object $M$ in $\sK$ together with a
morphism $M \otimes_\sK A \to M$ in $\BIMOD(\Ab,\sK)$ making the obvious
diagrams commute.
If $R$ is a commutative ring and $\sK=\Mod R$, then $R$-algebras in
the usual sense give $\sK$-algebras, and the new notion of a module
coincides with the old one.

Thus, the algebras in this paper are endo-functors with additional
structure, modules are objects in categories on which the functors
act, morphisms between algebras are natural transformations, et cetera.
Van den Bergh's language closely parallels the classical language of rings and modules
but one must be vigilant. Sometimes the proofs and statements of
results that are trivial for ordinary algebras and modules must be modified in subtle ways 
 to obtain the appropriate analogue for these ``functor
algebras''. For example, $A$ itself is not an $A$-module, but it
is an $A$-$A$-bimodule.

\subsection{Quasi-schemes \cite[Sect. 3.5]{vdB2}}

A {\sf quasi-scheme} $X$ is a Grothendieck category.
We will use the symbol $X$ when we think of the quasi-scheme as a geometric object
and $\Mod X$ when we think of it as a category. Objects in $\Mod X$
will be called $X$-modules.
We write $\mod X$ for the 
full subcategory of noetherian
$X$-modules. We say that $X$ is {\sf noetherian} if $\mod X$
generates $\Mod X$. Throughout $X$ will denote a noetherian 
quasi-scheme.

We use the word ``noetherian'' in the same way as Van den Bergh \cite[page 27]{vdB2} 
but it is not 
compatible with the usual notion of a noetherian scheme. If $X$ is a
noetherian scheme, it need not be the case that $\coh X$ 
generates $\Qcoh X$.

In this paper $X$ will play the role of the base space for a quantum projective
space bundle.

We write $o_X$ for the identity functor on $\Mod X$. 
It is an algebra in the sense of section \ref{sect.alg} and in some sense plays a
role similar to the structure sheaf or, more precisely, $- \otimes \cO_X$.

\subsection{Graded $X$-modules \cite[Sect 3.2]{vdB2}}
\label{sect.grmod}

A {\sf graded $X$-module} $M$ is a
sequence of objects $M=(M_n)_{n \in \ZZ}$ in $\Mod X$. We call
$M_n$ the degree $n$ component and, abusing notation, write $M=\oplus
M_n$. We define the category $\GrMod X$ by declaring a morphism
$f:M \to N$ to be a sequence of $X$-module maps $f_n:M_n \to N_n$.
We define $M_{>n}$ by declaring that $(M_{>n})_j=0$ if $j\le n$ and
$(M_{>n})_j=M_j$ if $j>n$.

\subsection{Graded $X$-algebras}
\label{sect.gralg}

{\sf Graded $X$-algebras}  $A=\oplus_{i \in\ZZ} A_i$ are defined  by combining the ideas in sections \ref{sect.alg} and \ref{sect.grmod} \cite[page 19]{vdB2}.
The  degree zero component $A_0$ of a graded algebra $A$ inherits an
$X$-algebra structure from $A$ and can also be given the
structure of a graded $X$-algebra concentrated in degree zero.
There is a graded $X$-algebra monomorphism $A_0 \to A$ and
 an augmentation map of graded $X$-algebras $A \to A_0$. The
composition $A_0 \to A \to A_0$ is the identity. 
A graded algebra $A$ is {\sf connected} if the structure map $o_X \to A_0$ is an isomorphism. 

A graded
$A$-module is a graded $X$-module with additional structure \cite[Sect. 3.2]{vdB2}. 
We write $\GrMod A$ for the
category of graded $A$-modules and $\grmod A$ for the full
subcategory of noetherian modules.

The identity functor on $\GrMod X$ can be given the
structure of a graded $X$-algebra. It is concentrated in degree
zero, and its degree zero component is isomorphic to $o_X$.
Thus, we write $o_X$ for both the $X$-algebra 
associated to $\Id_{\Mod X}$ and for the graded $X$-algebra
associated to $\Id_{\GrMod X}$.

\subsection{The quasi-scheme $\Projnc A$}
\label{sect.Proj}
Let $A$ be a connected graded $X$-algebra. We define $\Tors A$ to be the full subcategory of
 $\GrMod A$ consisting of  direct limits of graded $A$-modules that are annihilated by
  $A_{\ge n}$ for some $n \gg 0$. 
The full subcategory of noetherian objects in $\Tors A$ is denoted by $\tors A$.

 The {\sf projective quasi-scheme} associated to $A$ is denoted by $\Projnc A$
 and defined by declaring that 
$$
\Mod(\Projnc A) := \GrMod A/\Tors A.
$$
We write $\pi:\GrMod A \to \Projnc A$ for the quotient functor and $\omega$ for its right adjoint.

\subsection{$\cO_X$-bimodule algebras \cite{vdB1}}
\label{sect.OXY.bimods}

Let $X$ and $Y$ be noetherian schemes over a field $k$.
We write $\SSupp \cM$ for the scheme-theoretic support of a quasi-coherent
 $\cO_{X \times Y}$-module $\cM$.

A {\sf coherent $\cO _X$-$\cO _Y$ bimodule} is a coherent $\cO_{X \times Y}$-module $\cM$
 such that the projections
$pr _1:\SSupp \cM \to X$ and $pr _2:\SSupp \cM \to Y$ are finite.  A morphism of coherent
$\cO _X$-$\cO _Y$ bimodules is a morphism of
$\cO _{X\times Y}$-modules. 
  The category of  coherent $\cO _X$-$\cO _Y$
bimodules is denoted by $\Bimod(\cO_X,\cO_Y)$.

An {\sf $\cO_X$-bimodule algebra} is an algebra object in the
category of $\cO_X$-$\cO_X$-bimodules.
Graded $\cO_X$-bimodule algebras are then defined in an obvious way.
A graded $\cO_X$-bimodule algebra $\cB$ is analogous to a sheaf
of $\cO_X$-algebras, but it is actually a sheaf on $X \times X$
and local sections of $\cB$ do not have an algebra structure; its
global sections do, and they typically form a non-commutative graded
algebra.

\subsection{Quantum $\PP^n$-bundles (see section \ref{sect.Applic} for details)}
\label{sect.qu.Pn.bun}

The non-commutative analogue of a $\PP^n$-bundle over a scheme 
$X$ will be defined in the following way. First, take a coherent $\cO_{X \times X}$-module $\cE$ 
such that each $pr_{i*}\cE$ is locally free of finite rank, then use $\cE$ to construct a 
graded  $\cO_X$-bimodule algebra $\cB$ that 
is analogous to the symmetric algebra of a locally free $\cO_X$-module. 
The quasi-scheme $\PP(\cE)$ is defined
 implicitly by declaring that its module category is the quotient category 
$\Mod \PP(\cE) := \Projnc \cB:= \GrMod \cB/\Tors \cB$, where $\Tors \cB$
consists of the direct 
limits of graded modules that are annihilated by a power of the
augmentation ideal.
This definition is due to Van den Bergh \cite{Vn} (see also the thesis
of David Patrick \cite{Pat}).
If $\cE$ is a locally free $\cO_X$-module where $X$ is identified with the diagonal in $X \times X$, then
$\cB$ is the symmetric algebra $S(\cE)$ and $\Mod \PP(\cE)$
is $\Qcoh (\Proj \cB)$, the quasi-coherent sheaves on the
usual projective space bundle associated to $\cE$.

\subsection{Notation for Grothendieck groups}
\label{sect.Gr.gps}

Our goal is to compute the Grothendieck group of $\mod \PP(\cE)$, the category of noetherian modules
on the noetherian quasi-scheme $\PP(\cE)$.

We want to use notation that is the same as, or at least compatible with, the classical notation so we adopt the following conventions:
\begin{enumerate}
\item{}
We will write $K_0'(\sC)$ for the Grothendieck group of an abelian category $\sC$ with the following
exceptions.
\item{}
If $X$ is a noetherian scheme we will write $K_0'(X)$ for  $K_0'(\coh X)$ and
\item{}
$K_0(X)$ for
the Grothendieck group of the additive category of coherent locally free
$\cO_X$-modules. 
\item{}
When $X$ is a  separated, regular, noetherian scheme, the natural map \break
$K_0(X) \to K_0'(X)$ is an isomorphism so we identify these groups.
\item{}
If $X$ is a noetherian quasischeme we will write $K_0'(X)$ for $K_0'(\mod X)$. In particular, we do this when $X = \PP(\cE)$ and when $X=\Projnc A$.
\item{}
Because every noetherian $\PP(\cE)$-module has a 
finite resolution by modules that behave like locally free sheaves $\PP(\cE)$ is regular in
a suitable sense, and we therefore write $K_0(\PP(\cE))$ for $K_0'(\PP(\cE))$.
\item{}
When $X$ is a separated, regular, noetherian scheme and 
$A$ a flat, noetherian, connected graded $X$-algebra such that $\cTor^A_q(-,o_X)=0$ for $q \gg 0$
we will write $K_0(\Projnc A)$ for $K_0'(\Projnc A)$. The justification for doing this is explained
prior to Theorem \ref{thm.KY}.
\end{enumerate}
If $M$ is in $\sC$, we write $[M]$ for its class in $K_0'(\sC)$.

\subsection{Computation of $K_0(\Projnc A)$}
\label{sect.K0.proj}

We were unable to  compute $K_0'(\PP(\cE))$ by adapting the  method that Quillen used in \cite {Q} to compute the $K'_0$-group of a $\PP^n$-bundle.  Instead we follow Manin's method \cite {Man}. 

The first step towards that is to associate to the $\cO_X$-bimodule algebra $\cB$,
that  is  used to construct $\PP(\cE)$ (see section \ref{sect.qu.Pn.bun}), a  graded $X$-algebra $A$, 
in the sense of section \ref{sect.gralg}. We define $A$ by defining the left exact functors
$$
\cHom_X(A_i,-):=\cHom_{\cO_X}(\cB_i,-), \qquad i \ge 0.
$$
We show in section \ref{sect.Applic} that there is an equivalence of categories 
$$
\Mod \PP(\cE) \equiv \GrMod A/\Tors A 
$$
and therefore an isomorphism
$$
K_0'(\mod\PP(\cE)) \cong K_0'(\grmod A/\tors A)=K_0'(\Projnc A).
$$  
In section \ref{sect.K0} we compute $K_0'(\Projnc A)$ for a somewhat larger class of graded $X$-algebras. We allow $X$ to be any noetherian quasi-scheme, and $A$ any
flat, noetherian, connected, graded$X$-algebra such that $\ucTor^A_q(-,o_X)=0$ for $q \gg 0$  (i.e., 
the ``fibers'' of $\Projnc A \to X$ are ``smooth''). We show  that 
$$
K'_0(\grmod A) \cong K'_0(\grmod X)  \cong K'_0(X)[T,T^{-1}]
$$
and
$$
K'_0(\grmod A/\tors A) \cong K'_0(X)[T,T^{-1}]/I
$$
 where $I$ is the image of the map
$$
\rho:K'_0(\grmod X) \to K'_0(\grmod X), \qquad 
\rho[V]:=\sum (-1)^q [\ucTor^A_q(V,o_X)],
$$ 
where on the right hand side of the
equation, $V$ is viewed as a graded $A$-module via the
augmentation map $A \to o_X$. Up to this point in the argument the
structure sheaf of $X$ plays no role: $X$ is only assumed to be
 a quasi-scheme so does not have a structure sheaf.

\subsection{Computation of $K_0(\PP(\cE))$}

In section \ref{sect.Applic} we return to the main line of the argument by applying the 
computation of $K_0'(\Projnc A)$ in section \ref{sect.K0} to the special $A$ that is defined
in terms of $\cB$, the $\cO_X$-bimodule algebra used to coinstruct $\PP(\cE)$. In particular, we now
assume that $X$ is a separated, regular, noetherian scheme.

 We also write $K_0(\Projnc A)$ for $K_0'(\grmod A/\tors A)$ (the
remarks prior to Theorem \ref{thm.KY} justify this notational change).

A key step is to show that  the ideal $I$ mentioned in section \ref{sect.K0.proj} is principal: we show that
$$
\rho[V]=[V].\rho[\cO_X]
$$
 for all $V \in \grmod X$, where on the
right hand side we use the usual product in $K_0(X)[T,T^{-1}]$.
Therefore
$$
K_0'(\PP(\cE)) \cong K_0'(\Projnc A) \cong {{ K_0(X)[T,T^{-1}]}\over{(\rho [\cO_X])}}.
$$
The construction of $\cB$ is such that  its ``trivial module'', namely $\cO_X$, has a Koszul-like
resolution (Definition \ref{rs}). Now $\rho[\cO_X]$ can be read
off from this resolution, and is equal to the inverse of the
Hilbert series of $\cB$, $\sum[\cB_n]T^n \in K_0(X)[[T]]$ where
$\cB_n$ is the degree $n$ component of $\cB$, and $[\cB_n]$
denotes the class of $\cB_n$ in $K_0(X)$ when it is viewed as a
right $\cO_X$-module; to be precise, $[\cB_n]=[pr_{2*}\cB_n]$
where $pr_2:X \times X \to X$ is the projection onto the second
factor.

\subsection{Quantum ruled surfaces}

In section \ref{sect.int.thy} we further specialize to the case where $X$ is a smooth projective curve and $pr_{1*}\cE$ and $pr_{2*}\cE$ are locally free of rank two. We then call $\PP(\cE)$ a {\sf quantum
ruled surface}.  We show there is a well-behaved intersection
theory on $K_0(\PP (\cE))$.  
We define a fiber and a section as appropriate elements of $K_0(\PP
(\cE))$.  As an explicit application, we show that ``fibers do
not meet'' and ``a fiber and a section meet exactly
once'', as in the commutative case.

\section{Connected graded algebras over a quasi-scheme}
\label{sect.2}

From now on we will assume that $A$ is {\sf connected}.
 The unit and augmentation maps, $\eta:o_X
\to A$ and $\ve:A \to o_X$ are maps of graded $X$-algebras and
$\ve\eta=\Id_{o_X}$. Associated to $\eta$ and $\ve$ are various
functors \cite[(3.8), (3.9), (3.10)]{vdB2}. For example,
associated to $\eta$ is the forgetful functor
\begin{equation}
\label{eq.forget}
(-)_X:\GrMod A \to \GrMod X,
\end{equation}
its left adjoint $- \otimes_X A$,
and its right adjoint $\ucHom_X(A,-)$.
Associated to $\ve$ is the functor $(-)_A:\GrMod X \to
\GrMod A$, its left adjoint $-\otimes_A o_X$, and its right adjoint
$\ucHom_A(o_X,-)$.
The functors $(-)_X$ and $(-)_A$ are exact because they have both
left and right adjoints. Therefore their right adjoints
$\ucHom_X(A,-)$ and $\ucHom_A(o_X,-)$ preserve injectives.
We observe for later use that if $-\otimes_X A$ is exact, then
$(-)_X$ preserves injectives.
The functor $(-)_A$ is fully faithful so we often view $\GrMod X$ as
a full subcategory of $\GrMod A$.

The functors $\ucHom_A(o_X,-)$ and $-\otimes_A o_X$ associated to
the algebra map $A \to o_X$
 may be composed with the functor $(-)_A$ to obtain
an adjoint pair of functors sending $\GrMod A$ to $\GrMod A$, thus
endowing $o_X$ with the structure of a graded $A$-$A$-bimodule.
With this point of view, the appropriate notation for these two
functors is still  $\ucHom_A(o_X,-)$ and $-\otimes_A o_X$. The
relation between these functors is set out below, where on the left
hand side of each equality $o_X$ is viewed as an $A$-$A$-bimodule,
and on the right hand side it is viewed as a graded $X$-algebra:
\begin{align*}
- \otimes_A o_X & = (- \otimes_A o_X)_A,
\\
(-\otimes_A o_X)_X & = - \otimes_A o_X,
\\
\ucHom_A(o_X,-) & =\ucHom_A(o_X,-)_A,
\\
\ucHom_A(o_X,-)_X &= \ucHom_A(o_X,-).
\end{align*}

To illustrate how these functors are used we prove the following
basic result.

\begin{lemma}
\label{lem.izuru}
Let $A$ be a connected graded $X$-algebra. Let $M \in \GrMod A$,
and suppose that $M_i=0$ for $i<n$. Then the map $M=M \otimes_A A
\to M \otimes_A o_X$ of graded $A$-modules induces an isomorphism
$M_n \cong (M \otimes_A o_X)_n$ of $X$-modules. The kernel of the
epimorphism $M \to M \otimes_A o_X$ is contained in $M_{> n}.$
\end{lemma}
\begin{pf}
The $A$-module structure on $M$ is determined by a morphism $h:M
\otimes_X A \to M$, and $M\otimes _Ao_X$ is, by definition, the
coequalizer in $\GrMod A$ of

\begin{equation}
\label{eq.coeq}
\begin{CD}
M\otimes _XA @>{h}>> M \\
@V{1\otimes \ve}VV @. \\
M\otimes _Xo_X \cong M.
\end{CD}
\end{equation}
Let $\theta_i:\Mod X \to \GrMod X$ be defined by $\theta_i(V)_j$
equals $V$ if $i=j$ and $0$ if $i \ne j$.
If $M_{<n}=0$, then for every injective $I\in \Mod X$,
\begin{align*}
\Hom _X((M\otimes _XA)_n,I) & =\Hom_{\GrMod X}(M\otimes _XA,\theta _n(I)) \\
& =\Hom _{\GrMod X}(M, \ucHom _X(A,\theta _n(I))) \\
& =\prod _i\Hom _X(M_i, \ucHom _X(A,\theta _n(I))_i) \\
& =\prod _i\Hom _X(M_i, \cHom _X(A_{n-i},I)) \\
& = \Hom _X(M_n, \cHom _X(A_0,I)) \\
& =\Hom _X(M_n, I)
\end{align*}
so $(M\otimes _XA)_n \cong M_n$.
So the degree $n$ part of the diagram (\ref{eq.coeq})
consists of two identity maps, whence $(M\otimes _Ao_X)_n
\cong M_n$.
\end{pf}

We call $-\otimes _XA:\GrMod X \to \GrMod A$ the {\sf induction
functor}. Let $V \in \GrMod X$.
The structure map $h:V \otimes_X A \otimes_X A \to V
\otimes_X A$ making $V \otimes_X A$ a graded $A$-module is $V
\otimes \mu$, where $\mu: A \otimes_X A \to A$ is the
multiplication map.
A module $M\in \GrMod A$ is {\sf induced} if $M \cong V\otimes _XA$
for some $V\in \GrMod X$.

We say that $A$ is {\sf noetherian} if $X$ is noetherian and the induction functor
$-\otimes _XA:\GrMod X \to \GrMod A$ preserves noetherian modules
\cite[page 27]{vdB2}. If $A$ is noetherian and $M \in \grmod A$,
then each $M_i$ is a noetherian $X$-module.

\begin{lemma}
\label{1}
Let $A$ be an $\NN$-graded $X$-algebra.
\begin{enumerate}
\item{}
Every $M\in \GrMod A$ is a quotient of an induced module.
\item{}
If $X$ is noetherian, and $A$ is noetherian, then every $M \in
\grmod A$ is a quotient of an induced module $V \otimes_X A$ with
$V$ a noetherian graded $X$-module.
\end{enumerate}
\end{lemma}
\begin{pf}
(1)
Let $M \in \GrMod A$, and let $h:M \otimes_X A \to M$ be the
structure map. Forgetting the $A$-structure on $M$, we can form the
induced module $M \otimes_X A$.
Since $M$ is an $A$-module $h \circ
(M \otimes \mu)=h \circ (h \otimes A)$, from which it follows that
$h$ is a map of graded $A$-modules. To show $h$ is an
epimorphism in $\GrMod A$ it suffices to show it is an
epimorphism in $\GrMod X$. Suppose that $u:M \to N$ is a map of
graded $X$-modules such that $u\circ h=0$. Since $(M,h)$ is an
$A$-module, $h \circ (M \otimes \eta)=\id_M$; composing this
equality on the left with $u$ we see that $u=0$, whence $h$ is epic.

(2)
By (1) there is an epimorphism
$\oplus _i(M_i\otimes _XA)\to M$
where $M_i\in \grmod X$ is concentrated in degree $i$.
Since $M$ is noetherian, there is a finite direct
sum that still maps onto $M$, i.e.,
$$
(\oplus _{i=a}^{b}M_i)\otimes _XA\cong
\oplus_{i=a}^{b}(M_i\otimes _XA)\to M
$$
is epic for some $a\leq b$.
By an earlier remark, each $M_i$ is a  noetherian $X$-module, so
$\oplus _{i=a}^{b}M_i\in \grmod X$.
\end{pf}

We will prove a version of Nakayama's Lemma.  Recall that if $k$ is
a field, and $A$ is a connected $k$-algebra, then a finitely
generated graded $A$-module $M$ is zero if and only if $M
\otimes_A k=0$.

The {\sf augmentation ideal} of $A$ is $\fm:= \ker(\ve:A \to
o_X)$. Since $(-)_A: \GrMod X \to \GrMod A$ has a left and a right
adjoint, $\GrMod X$ is a closed subcategory \cite[page 20]{vdB2}
of $\GrMod A$, and hence $\fm$ is an ideal, not just a weak ideal,
in $A$ \cite[Section 3.4]{vdB2}. If $M
\in \GrMod A$, we write $M\fm$ for the image of the composition
$M\otimes _A\fm \to M \otimes _AA\cong M$.

\begin{lemma}
[Nakayama's Lemma]
\label{lem.Nak1}
Let $A$ be a connected $X$-algebra, and $M\in \grmod A$.
Then the following are equivalent:
\begin{enumerate}
\item{}
$M=0$;
\item{}
$M\fm =M$;
\item{}
$M\otimes _A o_X=0$.
\end{enumerate}
\end{lemma}
\begin{pf}
Clearly that (1) implies both (2) and (3).

The exact sequence $0\to \fm \to A\to o_X\to 0$
in $\BiGr (A,A)$ induces an exact sequence
$$M\otimes _A\fm \to M\otimes _AA\cong M\to M\otimes _A o_X\to 0$$
in $\GrMod A$, so $M \fm=\ker ( M\to M\otimes _A o_X)$.
Hence (2) and (3) are equivalent.
Because $M$ is noetherian, the chain of submodules $
\ldots \subset M_{\ge n} \subset M_{\ge n-1} \subset \ldots$ eventually
stops. Hence, if $M\neq 0$, there exists $n$ such that $M=M_{\geq n}$
and $M_n\neq 0$.  It then follows from Lemma \ref{lem.izuru} that
$M \otimes_A o_X \ne 0$. Therefore (3) implies (1).
\end{pf}

We now wish to consider $o_X$ as an $A$-$A$-bimodule and examine
$\ucTor^A_q(M,o_X)$ for $M \in \GrMod A$. The bifunctors $\cTor$
are defined in \cite[Section 3.1]{vdB2}. Because $M$ is a
$\ZZ$-$A$-bimodule the general definition in \cite{vdB2} admits a
small simplification.

\begin{lemma}
\label{lem.tor}
Let $\sD$ and $\sE$ be Grothendieck categories.
Let $\cM\in \sD$ and $\cN\in \BiMod (\sD ,\sE)$.
Then $\cTor^\sD_q(\cM,\cN)$ is in $\sE$, and is uniquely determined
by the requirement that there is a functorial isomorphism
$$
\Hom_\sE(\cTor^\sD_q(\cM,\cN),I) \cong \Ext^q_\sD(\cM,\cHom_\sE(\cN,I))
$$
for all injectives $I$ in $\sE$.
\end{lemma}
\begin{pf}
If $\cM \in \sD$ is viewed as a $\ZZ$-$\sD$-bimodule,
then $\cTor^\sD_q(\cM,\cN)$ is defined to be the $\ZZ$-$\sE$-bimodule
such that $\cHom_\sE(\cTor^\sD_q(\cM,\cN),I) =
\cExt^q_\sD(\cM,\cHom_\sE(\cN,I))$ for all injective $I$ \cite[page 12]{vdB2}.
However, if $\cF \in \cD$, then $\cExt^q_\sD(\cM,\cF) = \Ext^q_\sD(\cM,\cF)$, where
$\Ext^q_D(-,\cF)$ is the usual $q^{\th}$ right derived functor
of $\Hom_{\sD}(-,\cF)$. This is a consequence of \cite[3.1.2(4)]{vdB2}, and the fact
that the inclusion $\sD \to \BIMOD(\Ab,\sD)$ is exact and preserves injectives.

Let $\Inj \sE$ be the full subcategory of injectives in $\sE$.
There is a unique left exact functor
$\bar F:\sE \to \Ab$ extending the
functor $F:\Inj \sE \to \Ab$ defined by
$$F(I):=\Ext ^i_{\sD}(\cM,\cHom _{\sE}(\cN,I)).$$
We will show that $\bar F$ is representable.
This is sufficient to prove the lemma.

Since $\bar F$ is left exact, it is enough to show that it commutes
with products because then it will have a left adjoint and that
left adjoint evaluated at $\ZZ$ will give the representing object.
Let $\{W_j\}$ be a family of objects in $\sE$, and for each $j$
take an injective resolution $0\to W_j\to J_j^{\bullet}$.  Since the
product functor $\prod$ is left exact and exact on injectives,
and the product of injectives is injective,
$0\to \prod W_j\to \prod J_j^{\bullet}$ is an injective resolution of
$\prod W_j$. Therefore $\bar F(\prod W_j)=
\Ker (F(\prod J_j^0)\to F(\prod J_j^1))$.
Since $\cHom _{\sE}(\cN,-)$ and $\Ext ^i_{\sD}(\cM,-)$ both commute with
products, so does $F$.  Hence
\begin{align*}
\bar F(\prod W_j) & =\Ker (\prod F(J_j^0)\to \prod F(J_j^1))
\\
&
=\prod \Ker (F(J_j^0)\to F(J_j^1))
\\
&
=\prod \bar F(W_j).
\end{align*}
That is, $\bar F$ commutes with products.
\end{pf}

There is a graded version of this lemma that may be
applied to the present situation.
Thus, if $M \in \GrMod A$, then
$\ucTor^A_q(M,o_X) \in \GrMod A$ is determined by the requirement
that
\begin{equation}
\label{eq.Tor}
\Hom_{\GrMod A}(\ucTor_q^A(M,o_X),I) = \Ext^q_{\GrMod
A}(M,\ucHom_A(o_X,I))
\end{equation}
for all injectives $I \in \GrMod A$.

We may also view $o_X$ as an $A$-$o_X$-bimodule, and examine
$\ucTor^A_q(M,o_X)$ in that context. As remarked in \cite[page
18]{vdB2}, we have
$$
\ucTor^A_q(M,o_X) = \ucTor^A_q(M,o_X)_X,
$$
where on the left $o_X$ is an $A$-$o_X$-bimodule and on the right
it is an $A$-$A$-bimodule.

\smallskip

There is a notion of minimal resolution in $\grmod A$.

\begin{definition}
\label{defn.min.res}
A complex $\{P_{\bullet},d_n\}$ in $\grmod A$ is {\sf minimal} if
$\Im d_n\subseteq P_n\fm$ for all $n$.
\end{definition}

By \cite[3.1.5]{vdB2}, $\ucTor^A_\bullet(-,o_X)$ is a
$\d$-functor.

\begin{lemma} \label{lem.tor.minl.res'} Let $A$ be a connected
noetherian $X$-algebra, and $M\in \grmod A$. If induced modules
are acyclic for $\ucTor^A_\bullet(-,o_X)$, then
$\ucTor_i^A(M,o_X)_X \in \grmod X$.  In particular, if $M$ admits
a minimal resolution of the form $P_\bullet=V_\bullet \otimes_X A
\to M$, there is an isomorphism of graded $X$-modules
$$
V_i \cong \ucTor_i^A(M,o_X)_X.
$$
\end{lemma}

\begin{pf}
By Lemma \ref{1}, there is an epimorphism $V_0\otimes _XA\to M$
for some $V_0\in \grmod X$.
Repeatedly applying this process to the kernels of the
epimorphisms, we can construct a noetherian induced resolution
$P_{\bullet}=V_{\bullet}\otimes _XA\to M$.  Since $P_i\otimes_A
o_X=V_i\otimes _XA\otimes _Ao_X\cong V_i\in \grmod X$, by the
acyclicity hypothesis, $\ucTor_i^A(M,o_X)_X= H_i(P_\bullet
\otimes_A o_X)=H_i(V_{\bullet})\in \grmod X$. If the resolution
$P_\bullet=V_\bullet \otimes_X A \to M$ is minimal, then the
differentials in $P_\bullet \otimes_A o_X\cong V_{\bullet}$ are
all zero, hence
$
\ucTor_i^A(M,o_X)_X\cong V_i.
$
\end{pf}

{\bf Warning.}
Graded modules over a connected algebra in our sense need
not have minimal
resolutions. For example, if $A$ is the polynomial ring $k[x,y]$
with $\deg x=0$ and $\deg y=1$, then $A$ is a connected
$\AA^1$-algebra, but $A/(x)$ does not have a minimal resolution in
the sense just defined. If $A_0$ is semi-simple, then every
left bounded graded $A$-module has a minimal resolution.

For our applications to a quantum projective space bundle
$\PP(\cE)\to X$ over a (commutative!) scheme $X$,
we will only need that one very special module has a minimal
resolution, namely $\cO_X$, and that is guaranteed by the very
definition of $\PP(\cE)$.

\section{Flat connected graded $X$-algebras}

We retain the notation from the previous section.
Thus $X$ is a noetherian quasi-scheme and
$A$ is a noetherian, connected, graded $X$-algebra.

\begin{definition}
We say that $A$ is a {\sf flat $X$-algebra} if
$-\otimes _XA_i$ is exact for all $i$.
\end{definition}

Throughout this section we assume that $A$ is a flat $X$-algebra.

Thus $A$ is analogous to a locally free sheaf of graded
algebras on a scheme $X$.

Lemma \ref{lem.acyclic} shows that the flatness of $A$
implies that induced modules are acyclic for $\ucTor ^A_\bullet(-,o_X)$.
This is used to characterize those $M \in \grmod A$ for which
$\ucTor^A_1(M,o_X)=0$.

\begin{lemma} \label{lem.ext}
Suppose that $A$ is a flat, noetherian, connected, graded $X$-algebra.
If $V\in \GrMod X$ and $M\in \GrMod A$,
then there is a functorial isomorphism
$$\Ext ^q_{\GrMod A}(V\otimes _XA,M)\cong \Ext ^q_{\GrMod X}(V,M_X).$$
\end{lemma}
\begin{pf}
If $0\to M\to I^{\bullet}$ is an injective resolution in
$\GrMod A$, then $0\to M_X\to I_X^{\bullet}$ is an injective
resolution in $\GrMod X$. Thus
\begin{align*}
\Ext ^q_{\GrMod A}(V\otimes _XA,M) & =
H^q(\Hom _{\GrMod A}(V\otimes _XA,I^{\bullet})) \\
& = H^q(\Hom _{\GrMod X}(V,I_X^{\bullet})) \\
& =\Ext ^q_{\GrMod X}(V,M_X).
\end{align*}
\end{pf}

\begin{lemma}
\label{lem.acyclic}
Suppose that $A$ is a flat, noetherian, connected, graded $X$-algebra.
Then induced modules are acyclic
for $\ucTor ^A_\bullet(-,o_X)$.
\end{lemma}
\begin{pf}
Let $I\in \GrMod A$ be injective.
By the definition of $\ucTor$,
$$
\Hom _{\GrMod A}(\ucTor _q^A(V\otimes _XA,o_X),I) =\Ext ^q_{\GrMod
A}(V\otimes _XA,\ucHom _A(o_X,I)).
$$
By Lemma \ref{lem.ext},  this is isomorphic to
$\Ext ^q_{\GrMod X}(V,\ucHom _A(o_X,I)_X)$.
In the proof so far $o_X$ has been viewed as a graded
$A$-$A$-bimodule. But $\ucHom _A(o_X,I)_X = \ucHom_A(o_X,I)$, where
on the right hand side $o_X$ is viewed as a graded $X$-algebra.
Since $ \ucHom_A(o_X,-)$ is right adjoint to the exact functor
$(-)_A$ it preserves injectives. Hence
$\ucHom_A(o_X,I)$ is injective and, if $q \ge 1$,
$\ucTor ^A_q(V\otimes _XA,o_X)=0$.
\end{pf}

If $k$ is a field, $A$ a connected $k$-algebra, and $M$ a finitely
generated graded $A$-module, then $\Tor_1^A(M,k)=0$ if and only if
$M$ is a free $A$-module. The next result is an analogue of this.

\smallskip
{\bf Terminology:}
Given a sequence of subobjects $\cdots \subset M_1 \subset M_2 \subset \cdots $ in an abelian
category, we call the subquotients $M_i/M_{i-1}$ the {\sf slices} of this sequence.

\begin{lemma}
\label{lem.Nak2}
Let $A$ be a flat, noetherian, connected, graded $X$-algebra.
For $M \in \grmod A$, $\ucTor^A_1(M,o_X)=0$
if and only if $M$ has a finite filtration by graded $A$-submodules
such that each slice is an induced module.
\end{lemma}
\begin{pf}
($\Leftarrow$)
By the long exact sequence for $\ucTor^A_\bullet(-,o_X)$, it suffices to
show that $\ucTor_1^A(V \otimes_X A,o_X)$ is zero.  But
this is true by Lemma \ref{lem.acyclic}.

($\Rightarrow$)
By Lemma \ref{1},
there are integers $a$ and $b$ for which the natural map
 $$
\oplus_{i=a}^b M_i\otimes_X A \to M
$$
is an epimorphism.  We will prove the result by induction on $b-a$.

Suppose that $b-a=0$. Then there are exact sequences
$$
0 \to K \to M_a \otimes_X A \to M \to 0
$$
and
$$
0=\ucTor^A_1(M,o_X) \to K \otimes_A o_X \to M_a \otimes_X A
\otimes_A o_X @>\varphi>>  M  \otimes_A o_X  \to 0
$$
of graded $X$-modules. But $M_a \otimes_X A \otimes_A o_X \cong M_a$
because $- \otimes_X A \otimes_A o_X$ is equivalent to the identity
functor,
and $ M  \otimes_A o_X \cong M_a$ by Lemma \ref{lem.izuru}, so
$\varphi$ is
monic. Thus $K \otimes_A o_X=0$. Since $M$ is noetherian, $M_a$
is a noetherian $X$-module. Since $A$ is noetherian, $M_a \otimes_X
A$ is a noetherian $A$-module, whence $K$ is a noetherian
$A$-module.
Hence by Nakayama's Lemma, $K=0$. Thus $M \cong M_a \otimes_X
A$ is induced.

Now suppose that $b-a \ge 1$. Write $C:=\coker(M_a \otimes_X A \to
M)$ and $K:=\im(M_a \otimes_X A \to M)$. Hence there are exact
sequences
\begin{equation}
\label{eq.KMC}
0 \to K \to M \to C \to 0
\end{equation}
and
$$
0=\ucTor_1^A(M,o_X) \to \ucTor_1^A(C,o_X) \to K \otimes_A o_X \to M
\otimes_A o_X \to C \otimes_A o_X \to 0.
$$
However, $M_a \otimes_X A \otimes_A o_X \to M \otimes_A o_X$ is
monic,
so $K \otimes_A o_X \to M \otimes_A o_X$ is
monic. Therefore $ \ucTor_1^A(C,o_X)=0$.

The top and right arrows in the diagram
$$
\begin{CD}
\bigoplus_{i=a}^b M_i \otimes_X A @>>> M
\\
@VVV @VVV
\\
\bigoplus_{i=a}^b C_i \otimes_X A @>>> C
\end{CD}
$$
are epimorphisms, so the bottom arrow is epic too.
However $C_a=0$, so the induction hypothesis implies
that $C$ has a finite filtration with each slice an induced
module.

Since $C$ has a filtration by induced modules,
$\ucTor^A_2(C,o_X)=0$. Therefore $\ucTor_1^A(K,o_X)=0$.
Since the morphism $K_a \otimes_X A \to K$ is epic the proof of the
case $b=a$ shows that $K$ is induced. It follows
from the exact sequence (\ref{eq.KMC}) that $M$ has a filtration of
the required form.
\end{pf}

\section{Grothendieck groups}
\label{sect.K0}

We retain the notation from the previous section. Thus $X$ is a noetherian
quasi-scheme and  $A$ is a flat, noetherian, connected, graded $X$-algebra.

The degree shift functors $M \mapsto M(-1)$ induce automorphisms of $K'_0(\grmod X)$
and $K_0'(\grmod A)$. 
Let $\ZZ[T,T^{-1}]$ be the ring of Laurent polynomials and make $K'_0(\grmod X)$  and 
$K'_0(\grmod A)$ into $\ZZ[T,T^{-1}]$-modules by defining
$$
[M].T:=[M(-1)].
$$

\begin{proposition}
{\rm 
(cf. \cite[Theorem 6, page 110]{Q})}
\label{prop.K.grA}
Let $A$ be a flat, noetherian, connected, graded $X$-algebra.
Suppose  $\ucTor_q^A(M,o_X)=0$ for $q \gg 0$  for all $M \in \grmod A$.
Then there is an isomorphism
$$
\begin{CD}
\theta: K'_0(\grmod X) @>{\sim}>> K'_0(\grmod A),  \quad \theta[V]:= [V \otimes_X A],
\end{CD}
$$
of $\ZZ[T,T^{-1}]$-modules whose inverse is 
 the map 
 \begin{equation}
\label{eq.rho}
\begin{CD}
\rho: K'_0(\grmod A) @>{\sim}>> K'_0(\grmod X),  \;
\rho [M]: = \sum (-1)^i [\ucTor^A_i(M,o_X)].
\end{CD}
\end{equation}
 \end{proposition}
\begin{pf}
Since $ -\otimes_X A$ is exact it induces a group homomorphism
$\theta$. It is clear that $\theta$ is a $\ZZ[T,T^{-1}]$-module homomorphism. 

To see that $\rho$ is well-defined first observe that
$\ucTor^A_i(M,o_X) \in \grmod X$ by Lemmas \ref{lem.tor.minl.res'}
and \ref{lem.acyclic}, and the sum (\ref{eq.rho}) is finite by
hypothesis. The long exact sequence for $\ucTor^A_\bullet(-,o_X)$ then
ensures that $\rho$ is defined on $K'_0(\grmod A)$.

Now we show that $\theta$ is surjective.
Let $M \in \grmod A$. By Lemma
\ref{lem.Nak2}, there is an epimorphism $L:=V \otimes_X A \to M$ for
some $V \in \grmod X$.  By Lemma \ref{lem.acyclic},
$\ucTor^A_j(L,o_X)=0$ for all $j \ge 1$.
By induction, there is a resolution of graded $A$-modules
$$
\cdots \to L_n \to \cdots \to L_0 \to M \to 0
$$
with each $L_i$ induced from a noetherian graded $X$-module.
Consider the truncation
$$
0 \to N \to L_{d-1} \to  \cdots \to L_0 \to M \to 0.
$$
Since $\ucTor^A_i(L_j,o_X)=0$ for $0 \le j \le d-1$ and $i \ge 1$,
dimension shifting gives
$$
\ucTor_1^A(N,o_X) \cong \ucTor^A_{d+1}(M,o_X)=0.
$$
By Lemma \ref{lem.Nak2}, $N$ has a filtration with slices being
induced modules. Each $L_j$ is an induced module so, in
$K'_0(\grmod A)$, $[M]$ is a linear combination of induced modules.
This proves the surjectivity.

If $V \in \grmod X$, then Lemma \ref{lem.acyclic} shows that
$$
(\rho\circ \theta)[V]=\sum (-1)^i[\ucTor^A_i(V \otimes_X A, o_X)] =
[V \otimes_X A \otimes_A o_X]=[V],
$$
so $\rho\theta=\id_{K'_0(\grmod X)}$. Therefore $\theta$ is
injective and $\rho$ is its inverse.
\end{pf}

\begin{remark}
\label{rem.kgx}
If $M$ is a noetherian graded $X$-module, then its degree $i$ homogeneous component, $M_i$, is
zero for $i \ll 0$ and for $i\gg 0$. Hence  there is an isomorphism of $\ZZ[T,T^{-1}]$-modules
\begin{align*}
K'_0(\grmod X)   \overset{\sim}{\longrightarrow}
\; & K'_0(X)[T,T^{-1}] := K'_0(X)
\otimes_{\ZZ} \ZZ[T,T^{-1}].
\\
[M]\longmapsto & \sum _i[M_i]T^i.
\end{align*}
From now on we will make the identification 
$$
K'_0(\grmod X) \equiv K'_0(X)[T,T^{-1}] 
$$
without further comment.
\end{remark}

\begin{corollary}
Under the hypotheses of Proposition \ref{prop.K.grA},
$$
K'_0(\grmod A) \cong K'_0(X)[T,T^{-1}].
$$
\end{corollary}
\begin{pf}
Combine Proposition \ref{prop.K.grA} and Remark \ref{rem.kgx}.
\end{pf}

We now turn to the computation of
$$
K'_0(\Projnc A) = K_0'(\grmod A /\tors A).
$$
The degree shift
functor on graded $A$-modules preserves $\Tors A$, so induces an
auto-equivalence of $\Projnc A$ that sends noetherian objects to
noetherian objects, and hence induces an automorphism of
$K'_0(\Projnc A)$. We make $K'_0(\Projnc A)$ a $\ZZ[T,T^{-1}]$-module by
defining
$$
[\cM].T:=[\cM(-1)].
$$

A map $f$ of quasi-schemes is an adjoint pair of functors
$(f^*,f_*)$ between their module categories \cite[Section
3.5]{vdB2}.
We define a map 
$$
f:\Projnc A \to X
$$
 of quasi-schemes as follows.
If $\cM \in \Projnc A$, we define $f_*\cM=(\omega\cM)_0$, the
degree zero component. If $V \in \Mod X$, we define
$f^*V=\pi(V \otimes_X A)$.

If $A$ has the properties specified in sections 5 and 6
we think of $\Projnc A$ as a non-commutative analogue of a
projective space bundle over $X$ with structure map $f:\Projnc A \to X$.

\begin{theorem}
\label{thm.KP}
Retain the hypotheses of Proposition \ref{prop.K.grA}.
The group homomorphism 
$$
K'_0(X)[T,T^{-1}] \to K'_0(\Projnc A), \qquad [\cF]T^i \mapsto [f^*\cF(-i)],
$$
induces a $\ZZ[T,T^{-1}]$-module isomorphism
$$
{{K'_0(X)[T,T^{-1}]} \over{I}} \cong K'_0(\Projnc A),
$$
where  $I$ is the image in $K'_0(\grmod X) \cong K'_0(X)[T,T^{-1}]$
of the restriction to $K'_0(\tors A)$ of the map $\rho$ in (\ref{eq.rho}).
\end{theorem}
\begin{pf}
The localization sequence for $K$-theory gives an exact sequence
$$
\begin{CD}
K'_0(\tors A) @>{\psi}>> K'_0(\grmod A) @>>> K'_0(\Projnc A) @>>> 0.
\end{CD}
$$
Using the action of the degree shift functor, this is a sequence
of $\ZZ[T,T^{-1}]$-modules. If $M \in \tors A$, then $M \fm^i=0$
for $i \gg 0$, so objects of $\tors A$ have finite filtrations
with slices belonging to $\grmod X$. By D\'evissage, the inclusion
$\grmod X \to \tors A$ induces an isomorphism $K'_0(\grmod X) \to
K'_0(\tors A)$ which gives the left-most vertical isomorphisms in the diagram
$$
\begin{CD}
 K'_0(\tors A) @>{\psi}>> K'_0(\grmod A) @>>> K'_0(\Projnc A) @>>> 0
\\
@A{\cong}AA @V{\cong}V{\rho}V
\\
K'_0(\grmod X) @.  K'_0(\grmod X)
\\
@A{\cong}AA @VV{\cong}V
\\
K'_0(X)[T,T^{-1}] @. K'_0(X)[T,T^{-1}],
\end{CD}
$$
But $\rho$ is an isomorphism, so
$K'_0(\Projnc A) \cong \coker(\rho \circ \psi) \cong K'_0(X)[T,T^{-1}]/ I$.
\end{pf}

For the rest of this section, we suppose that $X$ is
a noetherian scheme in the usual sense.

As usual, $\cO_X$ denotes its structure sheaf.
We define $\Mod X$ to be the category of quasi-coherent
$\cO_X$-modules. Thus $\cO_X$ is in $\grmod X$, and hence in
$\grmod A$ via the map $A \to o_X$ of graded $X$-algebras and the
associated fully faithful functor $(-)_A:\GrMod X \to \GrMod A$.
Thus $\cO_X$ is a torsion $A$-module concentrated in degree zero.

\begin{proposition}
\label{prop.mod.map}
Retain the hypotheses of Proposition \ref{prop.K.grA} and let $\rho$ be the map defined by 
(\ref{eq.rho}).
Further, suppose that $X$ is a separated regular noetherian scheme.
Let $\cF \in \grmod X$, and view $\cF_A=\cF$ as a graded $A$-module
annihilated by $A_{> 0}$. Then
$$
\rho[\cF_A]=[\cF].\rho[\cO_X]
$$
where the product on the right hand side is the usual product in
$K'_0(X)$ extended to $K'_0(X)[T,T^{-1}]$ in the natural way.
\end{proposition}
\begin{pf}
The hypothesis on $X$ ensures that $\cF$ has a finite resolution by
locally free $\cO_X$-modules, say $\cW_\bullet \to \cF$.
The ring structure on $K'_0(X)[T,T^{-1}]$ is defined by
$$
[\cF].[\cG]= \sum (-1)^i [\uTor^X_i(\cF,\cG)]
$$
where $\uTor^X$ denotes the sheaf Tor of graded $X$-modules.

By definition of $\rho[\cO_X]$, we have
\begin{align*}
[\cF].\rho[\cO_X] &= \sum_q (-1)^q [\cF].[\ucTor^A_q(\cO_X,o_X)]
\\
&= \sum_q (-1)^q \sum_p (-1)^p
[\uTor^X_p(\cF,\ucTor^A_q(\cO_X,o_X))],
\end{align*}
Hence, to prove the proposition, we must show that
$$
\sum_n (-1)^n [\ucTor^A_n(\cF,o_X)] =  \sum_{p,q} (-1)^{p+q}
[\uTor^X_p(\cF,\ucTor^A_q(\cO_X,o_X))].
$$
This will follow in the standard way if we can prove the existence
of a convergent spectral sequence
\begin{equation}
\label{eq.spec.seq}
\uTor^X_p(\cF,\ucTor^A_q(\cO_X,o_X)) \Rightarrow \ucTor^A_{p+q}(\cF,o_X).
\end{equation}
This is a mild variation of the usual spectral sequence---the
details follow.

By Lemma \ref{1}, there is a resolution of $\cO_X$ in $\grmod A$
by noetherian induced modules.  Since every $\cO_X$-module is a
quotient of a locally free $\cO_X$-module there is a resolution
$\cV_\bullet \otimes_X A \to \cO_X$ with each $\cV_q$ a
locally free graded $\cO_X$-module. The spectral sequence we want will
come from the bicomplex $\cW_\bullet \otimes_X \cV_\bullet$.

First we make a preliminary calculation.
Applying the forgetful functor $(-)_X$, we can view $\cV_\bullet
\otimes_X A$ as an exact sequence in $\GrMod X$. Since $- \otimes_X
A$ is an exact functor, each $\cV_q \otimes_X A$ is a locally free
$\cO_X$-module. Hence $\uTor^X_q(\cF,\cO_X)\cong H_q(\cF \otimes_X \cV_\bullet
\otimes_X A)$. But $\uTor^X_q(\cF,\cO_X)=\cTor^X_q(\cF,\cO_X)$  is
zero if $q \ne 0$, so we conclude
that $\cF \otimes_X \cV_\bullet \otimes_X A \to \cF$ is exact, whence it
is a resolution of $\cF$ in $\grmod A$ by induced modules. Applying
$-\otimes_A o_X$ to this complex, it follows from
Lemma \ref{lem.acyclic} that
$\ucTor^A_q(\cF,o_X) \cong H_q(\cF \otimes_X \cV_\bullet)$.
Similarly, $H_q(\cV_\bullet)\cong \ucTor^A_q(\cO_X,o_X)$.

Now we consider the bicomplex  $\cW_\bullet \otimes_X \cV_\bullet$.
Since $-\otimes_X \cV_q$ is exact, taking homology along the rows
gives zero everywhere except in the first column where the homology
is $\cF \otimes_X \cV_\bullet$; taking homology down this column gives
$\ucTor_n^A(\cF,o_X)$ by the previous paragraph.
On the other hand, because $\cW_p$ is locally free, when we take
homology down the $p^{\th}$ column we get
$\cW_p \otimes_X H_\bullet(\cV_\bullet) \cong \cW_p \otimes_X
\ucTor^A_\bullet(\cO_X,o_X)$. If we now take the homology of
$ \cW_\bullet \otimes_X\ucTor^A_\bullet(\cO_X,o_X)$ we get
$\uTor^X_p(\cF,\ucTor^A_q(\cO_X,o_X))$. Thus, we obtain the desired
spectral sequence.
\end{pf}

\begin{remark}
\label{rem.spec.seq}
A connected graded algebra over a field $k$ has
finite global dimension if and only if the trivial module $k$ has
finite projective dimension; equivalently, if and only if
$\Tor_q^A(k,k)=0$ for $q \gg 0$. The spectral sequence
(\ref{eq.spec.seq}) yields an analogue of this result.
If $X$ is a separated regular noetherian scheme,
then $\ucTor^A_n(\cF,o_X)=0$ for all $\cF \in \mod X$ and all $n \gg 0$
if and only if $\ucTor^A_q(\cO_X,o_X)=0$ for $q \gg 0$.
\end{remark}

The next result is obtained by combining Theorem \ref{thm.KP} and
Proposition \ref{prop.mod.map}. First, however, we make a change of notation.

When $X$ is a separated, regular, noetherian scheme, every
coherent $\cO_X$-module has a finite resolution by locally free
$\cO_X$-modules, so the natural map $K_0(X) \to K_0'(X)$ from the
Grothendieck group of locally free coherent $\cO_X$-modules is an
isomorphism. When $X$ has these properties, and $A$ is a flat,
noetherian, connected, graded $X$-algebra such that
$\ucTor_q^A(-,o_X)=0$ for $q \gg 0$, every noetherian graded
$A$-module has a finite resolution by iterated extensions of
modules of the form $\cV \otimes_X A$ where $\cV$ is a locally
free coherent $\cO_X$-module (Proposition \ref{prop.K.grA}).
Since an extension by modules of the form $\cV \otimes_X A$ is
acyclic for $\ucTor^A_\bullet(-,o_X)$ it behaves like a flat
$A$-module. Although we do not have a notion of locally free for
$\Projnc A$, the image of a flat $A$-module in $\Projnc A$ 
behaves like a  locally free module. So  every noetherian module over
$\Projnc A$ has a finite resolution by modules that behave like locally free modules. We therefore write
$K_0(\Projnc A)$ for $K_0'(\Projnc A)$  in this situation.

\begin{theorem}
\label{thm.KY}
Let $X$ be a separated regular noetherian scheme and $A$ 
a flat, noetherian, connected, graded $X$-algebra in the sense of section \ref{sect.alg}.
Suppose that   $\ucTor_q^A(M,o_X)=0$ for $q \gg 0$  for all $M \in \grmod A$.
Then the group homomorphism 
$$
K_0(X)[T,T^{-1}] \longrightarrow K_0(\Projnc A), \qquad [\cF]T^i \longmapsto [f^*\cF(-i)],
$$
induces a $\ZZ[T,T^{-1}]$-module isomorphism
$$
K_0(\Projnc A) \cong {{K_0(X)[T,T^{-1}]}\over{(F)}}
$$
where $F$ is the polynomial
$\rho[\cO_X]=\sum_q(-1)^q[\ucTor_q^A(\cO_X,o_X)]$.
\end{theorem}
\begin{pf}
 With the identification $K_0'(\grmod X) \equiv K_0(X)[T,T^{-1}]$, the large diagram in the proof of 
 Theorem \ref{thm.KP} becomes
 $$
\begin{CD}
 K'_0(\tors A) @>{\psi}>> K'_0(\grmod A) @>>> K_0(\Projnc A) @>>> 0
\\
@V{\cong}VV @V{\cong}V{\rho}V 
\\
K_0(X)[T,T^{-1}]    @>>{\psiol} >  K_0(X)[T,T^{-1}]
\end{CD}
$$
and Proposition \ref{prop.mod.map} shows that the map $\psiol$ is multiplication by  $\rho[\cO_X]$.
The result follows. 
\end{pf}

{\bf Hilbert series.}
We say that $M \in \GrMod X$ is {\sf bounded below} if $M_i=0$ for $i \ll 0$. In particular, since $A$
is a noetherian graded $X$-algebra such that $A_i=0$ for $i<0$, every $M \in \grmod A$ is bounded
below.  Hence the forgetful functor 
$(-)_X:\grmod A \to \GrMod X$ induces a well-defined map
\begin{align*}
H:\grmod A & \longrightarrow K'_0(X)((T)),
\\
M & \longmapsto
H(M):=\sum _i[M _i]T^i
\end{align*}
 taking values in the the ring of formal Laurent series. 
Since $(-)_X$ is exact, $H$ induces a homomorphism
$$
K'_0(\grmod A) \longrightarrow K'_0(X)((T))
$$
of $\ZZ[T,T^{-1}]$-modules.
We call $H(M)$ the {\sf Hilbert series} of $M$.

Now suppose the hypotheses of  Theorem \ref{thm.KY} hold.
Evaluating $-\otimes_X A$ at $\cO_X$
gives $\cO_X \otimes_X A \in \grmod A$
with degree $i$ component $\cO_X \otimes_X A_i \in \mod X$, so
$$
H(\cO_X \otimes_A A)= \sum_{i=0}^\infty [\cO_X \otimes_X A_i]T^i \in  K_0(X)[[T]].
$$
Since $A$ is connected the leading coefficient of this is
$[\cO_X \otimes_X o_X]=[\cO_X]$ which is the identity element in the ring
$K_0(X)[[T]]$. Therefore $H(\cO_X \otimes_X A)$ is a unit
in $K_0(X)[[T]]$.

\begin{proposition}
\label{prop.Hseries.rho}
Under the hypotheses of  Theorem 
\ref{thm.KY}, if $\cO _X$ has a minimal resolution in $\grmod A$ by noetherian
induced modules, then
$$
\rho[\cO_X]=H(\cO_X \otimes_X A)^{-1}.
$$
\end{proposition}
\begin{pf}
Let $\cV_\bullet \otimes_X A \to \cO_X$ be a  minimal resolution in
$\grmod A$ (Definition \ref{defn.min.res}). It is necessarily finite.
Associativity of tensor product for weak bimodules
gives $\cV_q \otimes_X A=\cV_q \otimes_{\cO_X} (\cO_X \otimes_X A)$.
Therefore, we have the following computation in $K_0'(\grmod X) \subset K_0(X)((T))$:
\begin{align*}
[\cO_X]&=\sum(-1)^q[\cV_q \otimes_X A]
\\
&= \sum(-1)^q [\cV_q].[\cO_X \otimes A]
\\
&= \sum (-1)^q [\ucTor^A_q(\cO_X,o_X)].H(\cO_X \otimes A)
\\
&=\rho[\cO_X].H(\cO_X \otimes A).
\end{align*}
This proves the propoosition.
\end{pf}

Theorem \ref{thm.KY} and Proposition \ref{prop.Hseries.rho}
apply to ordinary algebras. If $k$ is a
field and $A$ is a connected graded noetherian $k$-algebra of
finite global dimension, then $A$ is an algebra over $X=\Spec k$
such that $\underline {\Tor} ^A_q(-,k)=0$ for $q\gg 0$, and
$-\otimes _kA_i$ is exact for every $i$.
Since $\rho [k]=H(A)^{-1}$,
where $H(A)$ is the usual Hilbert series of $A$,   we get
$$
K_0 (\Projnc A)\cong \ZZ [T]/(H(A)^{-1}),$$
thereby recovering \cite[Theorem 2.4]{MS}. In particular, this computes the $K_0$ for 
various non-commutative analogues of $\PP^n$.

\section{Applications to quantum projective space bundles over a commutative scheme}
\label{sect.Applic}

Now we apply the results of the previous section to quantum projective
space bundles over a commutative scheme.  

We use the definitions in section \ref{sect.OXY.bimods}.
Thus in this section  $X$, $Y$, and $Z$ are noetherian schemes over a field $k$.
Let $pr _1$ and $pr _2$ denote the canonical projections
$pr _1:X\times Y\to X$ and $pr _2:X\times Y\to Y$.
Let $\Mod X$ and $\Mod Y$ be the categories of quasi-coherent sheaves
on $X$ and $Y$ respectively.

 The tensor product of a coherent $\cO _X$-$\cO _Y$ bimodule $\cM$ and a
 coherent $\cO _Y$-$\cO _Z$ bimodule $\cN$  is the coherent 
$\cO _X$-$\cO _Z$ bimodule, defined by
\begin{equation}
\label{eq.C}
\cM \otimes _{\cO _Y}\cN :=pr _{13*}(pr _{12}^*\cM \otimes _{\cO _{X\times Y\times Z}}pr _{23}^*\cN ).
\end{equation}
We may tensor an $\cO _X$-$\cO _Y$ bimodule with an $\cO _X$-module or $\cO _Y$-module in the appropriate order.  For example, if $\cF \in \Mod X$ and $\cM $ is an $\cO _X$-$\cO _Y$ bimodule, then
\begin{equation}
\label{eq.spC}\cF \otimes _{\cO _X}\cM :=pr _{2*}(pr _{1}^*\cF \otimes _{\cO _{X\times Y}}\cM )
\end{equation}
is in $\Mod Y$.  This determines a functor $-\otimes _{\cO _X}\cM :\Mod X\to \Mod Y$.  Formula (\ref{eq.spC}) is a special case of (\ref{eq.C}) since $\cF$ is naturally an $\cO _{\Spec k}$-$\cO _X$ bimodule.

Let $\cM$ be a coherent $\cO _X$-$\cO _Y$ bimodule and let $\G =\SSupp \cM$ be its scheme-theoretic support.  Since $pr_2:\G \to Y$ is finite, $pr_{2*}:\Mod \G \to \Mod Y$ has a right adjoint $pr_2^!:\Mod Y\to \Mod \G$.  We define  $\cHom _{\cO _Y}(\cM ,-):\Mod Y\to \Mod X$ by
\begin{equation}
\label{eq.Hom}
 \cHom _{\cO _Y}(\cM ,\cG):= \cHom_{\cO_Y}(pr_{2*}\cM,\cG) = 
 pr_{1*}(\cHom _{\cO _{X\times Y}}(\cM,pr_2^!\cG))
\end{equation}
for $\cG \in \Mod Y$.

A special case of (\ref{eq.spC}) arises when $\cF=\cO_X$. This
amounts to viewing $\cM$ as a right $\cO_Y$-module. More precisely,
the direct image is an exact functor
$pr_{2*}:\Bimod(\cO_X,\cO_Y) \to \Mod Y$. Of course, $pr_{2*}$ need not be exact on $\Mod {X
\times Y}$, but its restriction to $\Bimod(\cO_X,\cO_Y)$ is because
if $0 \to \cF_1 \to \cF_2 \to \cF_3 \to 0$ is an exact sequence in
$\Bimod(\cO_X,\cO_Y)$, then $\cup_{i=1}^3 \Supp \cF_i$ is
affine over $Y$.

\begin{lemma}
The rule $\cM \longmapsto \cHom_{\cO_Y}(\cM,-)$ determines a fully faithful functor 
$\Bimod(\cO _X,\cO _Y) \longrightarrow \BiMod (X,Y)$ that is compatible with the tensor products.
\end{lemma}
\begin{pf}
Let $\cM$ be a coherent $\cO _X$-$\cO _Y$ bimodule. First we show that the
functor $-\otimes _{\cO _X}\cM:\Mod X\to \Mod Y$ defined by
(\ref{eq.spC}) is  left adjoint to the functor
$\cHom _{\cO _Y}(\cM ,-): \Mod Y\to \Mod X$ defined by (\ref{eq.Hom}).

Let $\G =\SSupp \cM$. Let $pr_{2*}$ denote the restriction to $\Mod
\Gamma$.  Then
\begin{align*}
\Hom _{\cO _Y}(\cF \otimes _{\cO _X}\cM,\cG) & =\Hom _{\cO _Y}(pr_{2*}(pr_1^*\cF \otimes _{\cO _{X\times Y}}\cM),\cG) \\
& \cong \Hom _{\cO _{\G}}(pr_1^*\cF \otimes _{\cO _{X\times Y}}\cM,pr_2^!\cG) \\
& \cong \Hom _{\cO _{\G}}(pr_1^*\cF ,\cHom _{\cO _{X\times Y}}(\cM,pr_2^!\cG)) \\
& \cong \Hom _{\cO _X}(\cF ,pr_{1*}\cHom _{\cO _{X\times Y}}(\cM,pr_2^!\cG)) \\
& =\Hom _{\cO _X}(\cF ,\cHom _{\cO _Y}(\cM ,\cG)),
\end{align*}
 which proves the adjointness claim.
 
We may view the left exact functor $\cHom _{\cO _Y}(\cM ,-)$ as an object in the 
opposite of the category of left exact functors; as such it is an
$X$-$Y$ bimodule in the sense of section \ref{sect.bimod} that we will denote by $\cM $.

 The tensor product (\ref {eq.C}) of bimodules is associative \cite [Proposition 2.5]{vdB1}, so is compatible with the 
tensor product defined in section \ref{sect.bimod}.  A morphism of $\cO _X$-$\cO _Y$ bimodules $\cM \to \cN$ induces natural trasformations
$$\cHom _{\cO _Y}(\cN,-)\to \cHom _{\cO _Y}(\cM ,-)$$
and
$$-\otimes _{\cO_ X}\cM\to -\otimes _{\cO _X}\cN.$$

The functor is fully faithful by \cite [Lemma 3.1.1]{Vn}.
\end{pf}

By symmetry $\cM \otimes _{\cO _Y}-:\Mod Y\to \Mod X$ has a right adjoint, so $\cM $ also determines a $Y$-$X$ bimodule.

\begin{definition} A coherent $\cO _X$-$\cO _Y$ bimodule $\cE $ is {\sf locally free} if $pr_{1*}\cE $ and $pr_{2*}\cE $ are locally free on $X$ and $Y$ respectively.  If $pr_{1*}\cE $ and $pr_{2*}\cE $ are locally free of the same rank $r$, then we say that $\cE $ is {\sf locally free of rank $r$}.
\end{definition}

Let $\cM $ be a coherent $\cO _X$-$\cO _Y$ bimodule.  Then $pr _{1*}\cM $ (resp. $pr _{2*}\cM $) is locally free if and only if the functor $-\otimes _{\cO _X}\cM $ (resp. $\cM \otimes _{\cO _Y}-$) is exact.  Let $\cE$ be an $\cO _X$-$\cO _Y$ bimodule and $\cF$ an $\cO _Y$-$\cO _Z$ bimodule.  If $\cE $ and $\cF$ are locally free, so is $\cE \otimes _{\cO _Y}\cF$.

\bigskip

From now on we restrict to the case $X=Y$. 
An $\cO _X$-$\cO _X$ bimodule will be called $\cO_X$-bimodule for short.
 Let $\Delta$
denote the diagonal in $X \times X$. Then $\cO_\Delta$ is an
$\cO_X$-bimodule in an obvious way. We call it the {\sf trivial}
$\cO_X$-bimodule. It is isomorphic to $\cO_X$ on each side, and
the $\cO_X$-action is central.

\begin{definition}
\label{defn.bimod.alg}
A {\sf graded $\cO _X$-bimodule algebra} is an $\cO _X$-bimodule
$\cB $ that is a direct sum $\cB =\oplus _i\cB _i$ of
$\cO_X$-bimodules equipped with maps of $\cO _X$-bimodules
$\cB _i\otimes_{\cO_X}\cB _j\to \cB _{i+j}$ (the multiplication) and
$\cO _\Delta\to \cB _0$ (the unit) such that the usual diagrams
commute.  We say that $\cB$ is {\sf connected} if $\cB$ is
$\NN$-graded and $\cB _0\cong \cO _\Delta$ as $\cO _X$-bimodules.
\end{definition}

\begin{definition}
Let $\cB$ be a graded $\cO_X$-bimodule algebra.
An $\cO _X$-module $\cM $ that is a direct sum $\cM =\oplus _i\cM _i$
of $\cO _X$-modules is called a {\sf graded right $\cB$-module}
if it is equipped with maps of $\cO _X$-modules
$\cM _i\otimes _{\cO _X}\cB _j\to \cM _{i+j}$ (the action)
satisfying the usual compatibilities.  The category of graded
right $\cB$-modules will be denoted by $\GrMod \cB$.

A graded right $\cB$-module is {\sf torsion} if
$\cM_n=0$ for $n\gg 0$.  The subcategory consisting of direct
limits of torsion modules is denoted by $\Tors \cB$.
The {\sf quantum projective scheme} associated to $\cB$ is
$$\Projnc \cB:=\GrMod \cB /\Tors \cB.$$
\end{definition}

We say that $\cB$ is {\sf noetherian} if $\GrMod \cB$ is a locally
noetherian Grothendieck category.

\begin{definition}
Let $\cB$ be a graded $\cO_X$-bimodule algebra. A {\sf graded
$\cB$-$\cO_X$ bimodule} is a direct sum of $\cO_X$-bimodules,
$\cM=\oplus \cM_i$, endowed with maps $\cB_i
\otimes_{\cO_X} \cM_j\to \cM_{i+j}$ of $\cO_X$-bimodules making the
usual diagrams commute.

If $\cB$ is connected then the augmentation map $\ve:\cB \to
\cO_\Delta$ makes $\cO_\Delta$ a graded $\cB$-$\cO_X$-bimodule.
\end{definition}

\begin{definition}
\label{rs}
A {\sf quantum $\PP ^n$-bundle} over a commutative scheme
$X$ is a quasi-scheme of the form $\Projnc \cB$ where $\cB$ is a
connected graded noetherian $\cO _X$-bimodule algebra such that
\begin{itemize}
\item{}
each $\cB _i$ is a locally free coherent $\cO _X$-bimodule;
\item{}
there is an exact sequence
\begin{equation}
\label {eq.bndl.0}
0\to \cB \otimes_{\cO_X} \cE' _{n+1} (-n-1)\to \cdots \to \cB
\otimes_{\cO_X} \cE'_1 (-1)\to \cB \otimes_{\cO_X} \cE'_0
\to \cO _\Delta\to 0
\end{equation}
of graded $\cB$-$\cO_X$-bimodules in which each $\cE'_ i$ is a
locally free $\cO _X$-bimodule of rank ${{n+1}\choose{i}}$
concentrated in degree zero;
\item{}
there is an exact sequence
\begin{equation}
\label {eq.bndl}
0\to \cE _{n+1}\otimes _{\cO _X}\cB (-n-1)\to \cdots \to \cE_1 \otimes _{\cO _X}\cB (-1)\to \cE _0\otimes _{\cO _X}\cB \to \cO _\Delta \to 0
\end{equation}
of graded $\cO_X$-$\cB$-bimodules in which each $\cE_ i$ is a locally free $\cO _X$-bimodule of rank ${{n+1}\choose{i}}$ concentrated in degree zero.
\end{itemize}
\end{definition}

Usual projective space bundles are quantum $\PP^n$-bundles:
a locally free sheaf  $\cE$ on a scheme $X$ may be viewed as a
sheaf on $\Delta$, the usual symmetric algebra $S(\cE)$ may play
the role of $\cB$, and there are exact sequences of the prescribed
form with $\cE_q=\cE_q'=\wedge^q\cE$, the $q^{\th}$ exterior power.

\smallskip

Let $\cB$ be a graded $\cO _X$-bimodule algebra.  If each $\cB _i$ is a coherent $\cO _X$-bimodule, then $\cB$ determines a graded $X$-algebra $A$ in the earlier sense by
$$\cHom _X(A_i,-):=\cHom _{\cO_X}(\cB _i,-)$$
each of which has a left adjoint $-\otimes _XA_i=-\otimes _{\cO
_X}\cB_i$.  If $\cB $ is connected, so is $A$. A graded right $\cB
$-module $\cM $ can be thought of as a graded $A$-module because
$$\cM _i\otimes _XA _j=\cM _i\otimes _{\cO _X}\cB _j\to \cM _{i+j}.$$
Conversely, a graded $A$ module $\cN$ can be thought of as a graded right $\cB$-module because
$$
\cN _i\otimes _{\cO _X}\cB _j=\cN _i\otimes _XA _j\to \cN _{i+j}.
$$
It follows that $\GrMod \cB \cong \GrMod A$.
Moreover, this equivalence preserves torsion objects, so
$$
\Projnc \cB \cong \Projnc A.
$$

\begin{lemma}
\label{lem.acyclic2}
Let $X$ be a scheme, and $A$ a flat, noetherian, connected, graded $X$-algebra.
If $\cE$ is a locally free $\cO_X$-bimodule and $M \in \GrMod A$, then
$A \otimes_X \cE$ is acyclic for $\ucTor^A_\bullet(M,-)$.
\end{lemma}
\begin{pf}
We can view $\cE$ as a graded $X$-$X$-bimodule in the earlier sense,
concentrated in degree zero.
Composing left exact functors gives a graded $A$-$X$-bimodule
$A \otimes_X \cE$.
By the graded version of Lemma \ref{lem.tor},
$\ucTor_q^A(M,A \otimes_X \cE)$ is the graded $A$-module determined
by the requirement that
$$
\Hom_{\GrMod A}(\ucTor_q^A(M,A \otimes_X \cE),I) = \Ext^q_{\GrMod
A}(M,\ucHom_{X}(A \otimes_X \cE,I))
$$
for all injectives $I \in \GrMod X$. However,
$\ucHom_{X}(A \otimes_X \cE,-)$ has an exact left adjoint $-
\otimes_X A \otimes_{\cO_X}\cE$, so
$\ucHom_{X}(A \otimes_X \cE,I)$ is
injective. Therefore $\ucTor_q^A(M,A \otimes_X \cE)=0$ if
$q \ge 1$.
\end{pf}

\begin{theorem}
\label {thm.Kqpb}
Let $X$ be a separated regular noetherian scheme.
If $\Projnc\cB$ is a quantum $\PP ^n$-bundle over $X$ as in Definition
\ref{rs}, then the Hilbert series of $pr_{2*}\cB=
\cO_X \otimes_{\cO_X}\cB$ is
$$
H(pr_{2*} \cB)=
\biggl(\sum _q(-1)^q[pr_{2*}\cE _q]T^q \biggr)^{-1},
$$
and
$$K_0(\Projnc\cB)\cong K_0(X)[T]/( H(pr_{2*}\cB)^{-1}).$$
\end{theorem}
\begin{pf}
Let $A$ be the $X$-algebra determined by $\cB$.
Then $\Projnc A = \Projnc \cB$.
The result will follow from Theorem \ref{thm.KY} and Proposition
\ref{prop.Hseries.rho}, but in order to
apply those results we must check that $A$
satisfies the hypotheses of Proposition \ref{prop.K.grA}.

First $A$ is connected and noetherian because $\cB$ is. Since each
$\cB_j$ is locally free, $-\otimes_X A$ is exact, and $A$ is a flat $X$-algebra.
The sequence
(\ref{eq.bndl.0}) gives an exact sequence $A \otimes_X
\cE'_{\bullet} \to o_X \to 0$ in $\BiGr(A,X)$. By Lemma
\ref{lem.acyclic2}, we can compute $\ucTor^A_q(M,o_X)$ as the
homology of $M \otimes_A A \otimes_X  \cE'_{\bullet}$. But
(\ref{eq.bndl.0}) is a bounded complex, so $\ucTor^A_q(M,o_X)=0$
for $q \gg 0$. Thus, the results in section 4 apply to this
particular $A$.

Apply $pr_{2*}$ to the sequence (\ref{eq.bndl}). This gives an
exact sequence of graded right $\cB$-modules. Passing to the
category of $A$-modules this gives a resolution of $\cO_X$ by
induced modules, namely
$$
0\to pr_{2*}\cE _{n+1}\otimes _XA (-n-1)\to \cdots \to pr_{2*}\cE
_0\otimes _XA \to \cO _X\to 0.
$$
It is a minimal resolution because each $pr_{2*}\cE_i$ is
concentrated in degree zero. By Lemma
\ref{lem.tor.minl.res'}, $\ucTor ^{A}_q(\cO _X,o _X)\cong
pr_{2*}\cE _q(-q)$ so the result is a consequence of
 Theorem \ref{thm.KY} and Proposition \ref{prop.Hseries.rho}.
\end{pf}

As remarked after Definition \ref{rs},
usual projective space bundles are quantum $\PP^n$-bundles, so
Theorem \ref{thm.Kqpb} computes the Grothendieck group of a usual
projective space bundle.

\section{Intersection theory on quantum ruled surfaces}
\label{sect.int.thy}

In this section, let $X$ be a smooth projective curve
over an algebraically closed field $k$.  We apply our
computation of the Grothendieck groups to a quantum ruled surface
(defined below) and then use this to
develop an intersection theory for quantum ruled surfaces.

First we recall Van den Bergh's definition of a quantum ruled surface
\cite{Vn} (see also \cite {Ng}, \cite {Pat}, \cite {vdB1}).
Let $\cE$ be a locally free $\cO _X$-bimodule. By \cite[Section
3]{Vn}, $\cE$  has both a right and a left adjoint,
$\cE ^*$ and ${}^*\cE$, each of which is a locally free $\cO
_X$-bimodule. An invertible $\cO _X$-subbimodule $\cQ
\subset \cE \otimes _{\cO _X}\cE$ is {\sf nondegenerate} if the
composition
$$\cE ^*\otimes _{\cO _X}\cQ \to \cE ^*\otimes _{\cO _X}\cE \otimes _{\cO _X}\cE \to
\cE$$ induced by the canonical map $\cE ^*\otimes _{\cO _X}\cE \to
\cO _{\Delta}$, and the composition
$$\cQ \otimes _{\cO _X}{}^*\cE \to \cE \otimes _{\cO _X}\cE \otimes _{\cO _X}{}^*\cE \to
\cE$$ induced by the canonical map $\cE \otimes _{\cO _X}{}^*\cE
\to \cO _{\Delta}$ are isomorphisms.

\begin{definition}
\cite{Vn}
A quantum ruled surface over $X$ is a quasi-scheme $\PP
(\cE)$ where
\begin{itemize}
\item{} $\cE$ is a locally free $\cO _X$-bimodule of rank 2,
\item{} $\cQ\subset \cE \otimes _{\cO _X}\cE$ is a nondegenerate
invertible $\cO _X$-subbimodule, \item{} $\cB =T(\cE)/(\cQ)$ where
$T(\cE)$ is the tensor algebra of $\cE$ over $X$, and \item{}
$\Mod \PP (\cE)=\Projnc \cB $.
\end{itemize}
We define $\cO_{\PP(\cE)}:=\pi( pr_{2*}\cB)=\pi(\cO_X \otimes_X \cB)$.
\end{definition}

By \cite[4.2.1]{Vn}, $\PP (\cE)$ is independent of the choice of a
nondegenerate $\cQ$.  In fact, $\cQ$ is not even needed
to define $\PP (\cE)$ (loc. cit.).

\begin{theorem}
\label {thm.grrs}
If $\PP (\cE)$ is a quantum ruled surface over
$X$, then
\begin{equation}
\label{eq.Krs} K_0(\PP (\cE))\cong
{{K_0(X)[T]}\over{(F)}}
\end{equation}
where $(F)$ is the ideal generated by $F=[\cO_X]-[pr_{2*}\cE]T+[pr_{2*}\cQ]T^2.$
\end{theorem}

\begin{pf}
By \cite[Theorem 1.2]{Vn}, $\cB$ is noetherian.
By \cite[Theorem 6.1.2]{Vn}, each $\cB _n$  is a
locally free coherent $\cO _X$-bimodules of rank $n+1$ for
$n\geq 0$, and there is an exact sequence
$$0\to \cQ \otimes _{\cO _X}\cB (-2)\to \cE \otimes _{\cO _X}\cB (-1)\to \cO _{\Delta}\otimes _{\cO _X}\cB \to \cO _{\Delta}\to 0$$
of graded $\cO _X$-$\cB$ bimodules.   By symmetry of the
definition of $\cB$, there is an exact sequence
$$0\to \cB \otimes _{\cO _X}\cQ (-2)\to \cB \otimes _{\cO _X}\cE (-1)\to \cB \otimes _{\cO _X}\cO _{\Delta} \to \cO _{\Delta}\to 0$$
of graded $\cB$-$\cO _X$ bimodules.  It follows that $\PP (\cE)$
is a quantum $\PP ^1$-bundle over $X$ in the sense of Definition
\ref {rs}, so Theorem \ref {thm.Kqpb} applies and gives the desired
result.
\end{pf}

The {\sf structure map }
$$
f:\PP (\cE)\to X
$$
for a quantum ruled surface
over $X$ is given by the adjoint pair of
functors
\begin{align*}
f_*:\Mod \PP (\cE) & \to \Mod X
\\
f_*\cM & = (\o \cM)_0,
\\
f^*:\Mod X & \to \Mod \PP (\cE)
\\
f^* \cF &= \pi (\cF \otimes _{\cO _X}\cB).
\end{align*}
We have $f^* \cO_X=\cO_{\PP(\cE)}$ and $f_*\cO_{\PP(\cE)} =
\cO_X$. The next lemma, which is a simple consequence of results of
Nyman and Van den Bergh, shows among
other things that $R^1f_*\cO_{\PP(\cE)}=0$.

Notice that $f^*$ is exact because
$-\otimes _{\cO _X}\cB :\Mod X\to \GrMod \cB $ and $\pi
:\GrMod \cB \to \Mod \PP (\cE)$ are exact.
Moreover, $f^*\cF \in \mod \PP
(\cE)$ for all $\cF\in \mod X$ by \cite [Lemma 2.17]{Ns}, and
$R^qf_*\cM\in \mod X$ for all $ \cM\in \mod \PP (\cE)$
and $q\geq 0$ by \cite [Corollary 3.3]{Nsf}.

The isomorphism in Theorem \ref{thm.grrs} is induced by the map
$K_0(X)[T,T^{-1}] \to K_0(\PP(\cE))$, $[\cF]T^i \mapsto
[f^*\cF(-i)]$. Since $\cQ$ is invertible, $[pr_{2*}\cQ]$ is a unit
in $K_0(X)[T,T^{-1}]$, and we obtain an isomorphism $K_0(\PP(\cE))
\cong K_0(X) \oplus K_0(X)T$.

\begin{lemma}
\label {lem.hdi}
Let $\cL $ be a locally free coherent $\cO_X$-module. Then
\begin{align*}
& f_*((f^*\cL )(n))\cong \cL \otimes _{\cO _X}\cB _n \text { for all } n\in \ZZ, \\
& R^1f_*((f^*\cL)(n))\cong \begin{cases} 0 & \text { if } n\geq
-1,
\\
\cL \otimes
_{\cO _X}(\cB _{-2-n})^* & \text { if } n\leq -2, \\
\end{cases} \\
& R^if_*=0 \text { for } i\geq 2.
\end{align*}
Furthermore, $f^*f_*=\id$, $R^1f_* \circ f^*=0$, and
$Rf_* \circ f^* =\id_{{\sf D}^b(X)}$, where ${\sf
D}^b(X)$ denotes the bounded derived category of coherent
$\cO_X$-modules.
\end{lemma}
\begin{pf}
By \cite [Theorem 2.5]{Nsf},
\begin{align*}
& \tau (\cL \otimes _{\cO _X}\cB )=0, \\
& R^1\tau (\cL \otimes _{\cO _X}\cB )=0, \\
& (R^2\tau (\cL \otimes _{\cO _X}\cB ))_n\cong \begin{cases} 0 &
\text { if } n\geq
-1,  \\
\cL \otimes _{\cO _X}(\cB _{-2-n})^* & \text { if } n\leq -2, \\
\end{cases} \\
& R^i\tau =0 \text { for } i\geq 3.
\end{align*}
If $\cM \in \GrMod \cB $, there is an exact sequence
$$0\to \tau \cM \to \cM  \to \omega \pi \cM  \to R^1\tau (\cM )\to 0,$$
and isomorphisms $R^i\omega (\pi \cM )\cong R^{i+1}\tau \cM $ for
$i\geq 1$ in $\GrMod \cB $ \cite [section 3.8]{vdB2}.  Since
$R^if_*(\pi \cM (n))\cong R^i\omega (\pi \cM )_n$ for $i\geq 0$,
the result follows.

Let $\cF$ be a coherent $\cO_X$-module. There is an exact sequence
$0 \to \cL' \to \cL \to \cF \to 0$ with $\cL'$ and $\cL$ locally
free $\cO_X$-modules.
Because $f^*$ is exact, $f_*f^*$ is left exact; because $R^1f_*
\circ f_*$ vanishes on locally free sheaves, there is an exact
sequence $0 \to f_*f^* \cL' \to f_*f^* \cL \to f_*f^* \cF \to 0$.
Because $f_*f^*$ is the identity on locally free sheaves it follows
that $f_*f^* \cF \cong \cF$. Hence $f_*f^*$ is isomorphic to the
identity functor. It follows that
$Rf_* \circ f^* =\id_{{\sf D}^b(X)}$.

Because $R^2f_*=0$, $R^1f_*$ is right exact. Because $R^1f_* \circ
f^*$ vanishes on locally free sheaves it vanishes on all
$\cO_X$-modules.
\end{pf}

Because $R^2f_*=0$, there is a group homomorphism
$Rf_*:K_0(\PP(\cE)) \to K_0(X)$.

\medskip

{\bf Intersection theory.}
Let $Y$ be a noetherian quasi-scheme over a field $k$ such that
\begin{itemize}
\item{} $\dim_k \Ext^i_Y(\cM, \cN) <\infty$ for all $i\geq 0$, and
\item{} $\Ext^i_Y(\cM, \cN)=0$ for all $i \gg 0$,
\end{itemize}
for all $\cM , \cN\in \mod Y$.
A quantum ruled surface satisfies the first of these conditions
by \cite [Corollary 3.6]{Nsf},
and satisfies the second by Proposition \ref{prop.hd} below.
We define the Euler form on the Grothendieck group of $Y$ by
$$
(-,-) :K_0(Y) \times K_0(Y) \to \ZZ.
$$
$$
([\cM],[\cN]):=\sum_{i=0}^\infty (-1)^i \dim_k \Ext^i_Y(\cM, \cN).
$$
We define the intersection multiplicity of $\cM$ and $\cN$ by
$$\cM
\cdot \cN :=(-1)^{\dim \cN}(\cM , \cN)$$
for some suitably defined dimension function $\dim \cN $. In
particular, if $\cM$ and $\cN$ are ``curves" on a ``quantum
surface" $Y$, then we define $\cM \cdot \cN =-(\cM , \cN)$.

The {\sf homological dimension} of $\cM \in \mod Y$ is $$\hd
(\cM)=\sup \{ i \; | \; \Ext ^i_Y(\cM , \cN)\neq 0 \text { for
some } \cN \in \mod Y\},$$ and the {\sf homological dimension} of $Y$
is
$$\hd (Y)=\sup \{\hd (\cM) \; | \; \cM \in \mod Y\}.$$

\begin{proposition}
\label {prop.hd}
The quantum ruled surface $\PP (\cE)$ over $X$ has finite homological dimension.
\end{proposition}
\begin{pf}
First we show that $\hd (f^*\cF)\leq 2$ if $\cF$ is a coherent $\cO_X$-module.
Let $\cM \in \mod \PP(\cE)$.

Since $f_*$ is right adjoint to the exact functor $f^*$, it
preserves injectives, so there is a Grothendieck spectral sequence
$$
\Ext_X^p(\cF, R^qf_* \cM ) \Rightarrow \Ext^{p+q}_{\PP
(\cE)}(f^*\cF, \cM )
$$
for any $\cF \in \mod X$ and $ \cM  \in \mod \PP (\cE)$.
 Since $R^qf_*=0$ for $q\geq 2$ by Lemma \ref{lem.hdi}, we get an
exact sequence
\begin{align*}
0 \to \Ext^1_X(\cF, f_* \cM ) & \to \Ext^1_{\PP(\cE)}(f^*\cF,
 \cM ) \to \Hom_X(\cF, R^1f_* \cM ) \to
\\
\to \Ext^2_X(\cF, f_* \cM ) & \to \Ext^2_{\PP (\cE)}(f^*\cF,
 \cM ) \to \Ext^1_X(\cF, R^1f_*\cM ) \to
\\
\to \Ext^3_X(\cF, f_* \cM ) & \to \Ext^3_{\PP (\cE)}(f^*\cF,
 \cM ) \to \Ext^2_X(\cF, R^1f_*\cM) \to \cdots .
\end{align*}
Since $X$ is a smooth curve, $\Ext^p_X=0$ for $p\geq
2$, so $\Ext ^p_{\PP (\cE)}(f^*\cF, -)=0$ for $p\geq 3$ and
$\hd (f^*\cF)\leq 2$.

By Theorem \ref{thm.grrs},
$K_0(\PP (\cE))$ is generated by $[(f^*\cF)(n)]$ for $\cF\in \mod
X$ and $n\in \ZZ$. Let $0\to  \cM'\to  \cM \to  \cM''\to
0$ be an exact sequence in $\mod \PP (\cE)$. If two of $\cM,
 \cM', \cM''$ have finite homological dimension, so
does the other, so it is enough to show that $(f^*\cF)(n)$ has
finite homological dimension. But we have just seen that
$\hd ((f^*\cF)(n))\leq 2$ so the result follows.
\end{pf}

\begin{lemma}
\label {lem.prf}
Let $\PP (\cE)$ be a quantum ruled surface over $X$. If $\cF\in \mod X$
and $\cM \in \mod \PP(\cE)$, then
\begin{enumerate}
\item{}
$(f^*\cF,  \cM ) =(\cF, f_* \cM )-(\cF, R^1f_*\cM)$ and
\item{}
if $\cL$ is a locally free coherent $\cO _X$-module and
$n\geq -1$, then
$$(f^*\cF, (f^*\cL)(n)) =(\cF, \cL \otimes _{\cO _X}\cB _n).$$
\end{enumerate}
\end{lemma}
\begin{pf}
(1)
From the long exact sequence in the previous proof, we obtain
an exact sequence
$$
0 \to \Ext^1_X(\cF, f_* \cM )  \to \Ext^1_{\PP (\cE)}(f^*\cF,
 \cM ) \to \Hom_X(\cF, R^1f_*\cM ) \to 0
$$
and an isomorphism
$$
\Ext^2_{\PP(\cE)}(f^*\cF,  \cM ) \cong \Ext^1_X(\cF, R^1f_*\cM).
$$
Combining this with the adjoint isomorphism
$$\Hom_{\PP (\cE)}(f^*\cF,  \cM ) \cong
\Hom_X(\cF, f_* \cM ),$$
we get 
$$(f^*\cF ,  \cM ) = (\cF, f_* \cM ) - (\cF, R^1f_*\cM).$$

(2)
By Lemma \ref {lem.hdi},
$f_*((f^*\cL )(n))\cong \cL \otimes _{\cO _X}\cB _n$ and
$R^1f_*((f^*\cL )(n))=0$ so part (1) gives
$
(f^*\cF, (f^*\cL )(n)) =(\cF, \cL \otimes _{\cO
_X}\cB _n).$
\end{pf}

Our remaining goal is Theorem \ref{thm.int.thy} which shows that the
intersection theory on $\PP (\cE)$ is like that for a commutative ruled
surface.

First we recall the following standard result that describes
the Euler form on the smooth projective curve $X$.

\begin{lemma} \label {lem.curve}
Let $g$ be the genus of $X$, and $p$ and $q$ closed points of $X$.
Then
\begin{align*}
& (\cO _X, \cO _X)=1-g. \\
& (\cO _X, \cO _p)=1. \\
& (\cO _p, \cO _X)=-1. \\
& (\cO _p, \cO _q)=0.
\end{align*}
\end{lemma}

\begin{lemma} \label {lem.ed}
Let $p$ be a closed point of $X$ and $n$ an integer $\ge -1$. Then
\begin{enumerate}
\item{}
$R^1f_*((f^*\cO _p)(n))=0$,
\item{}
$[f_*((f^*\cO _p)(n))]=\sum _{i=0}^n[\cO _{p_i}]$ in $K_0(X)$
for some closed points $p_0, \cdots , p_n\in X$, and
\item{}
$(\cO _{\PP(\cE)}, (f^*\cO _p)(n)) =n+1.$
\end{enumerate}
\end{lemma}
\begin{pf}
Let $n \in \ZZ$.
Since $f^*$ is exact, the exact sequence
$0\to \cO _X(-p)\to \cO _X\to \cO _p\to 0$
in $\mod X$ induces an exact sequence
$$0\to (f^*(\cO _X(-p)))(n)\to (f^*\cO _X)(n)\to (f^*\cO _p)(n)\to 0$$
in $\mod \PP (\cE)$. Applying $f_*$ to this produces
a long exact sequence
\begin{align*}
0 & \to f_*((f^*\cO _X(-p))(n))\to f_*((f^*\cO _X)(n))\to f_*((f^*\cO _p)(n)) \\
& \to R^1f_*((f^*\cO _X(-p))(n))\to R^1f_*((f^*\cO _X)(n))\to
R^1f_*((f^*\cO _p)(n))\to \cdots
\end{align*}
in $\mod X$.

Now suppose that $n \ge -1$.

Applying Lemma \ref{lem.hdi} to various terms in the previous long
exact sequence, we get an exact sequence
$$0\to \cO _X(-p)\otimes _{\cO _X}\cB _n\to
\cO _X\otimes _{\cO _X}\cB _n\to f_*((f^*\cO _p)(n))\to 0,
$$
and also see that $R^1f_*((f^*\cO _p)(n))=0.$
On the other hand,
there is an exact sequence
$$0\to \cO _X(-p)\otimes _{\cO _X}\cB _n\to \cO _X\otimes _{\cO _X}\cB _n\to
\cO _p\otimes _{\cO _X}\cB _n\to 0$$ in $\mod X$, so
$$[f_*((f^*\cO _p)(n))]=[\cO _X\otimes _{\cO _X}\cB _n]-[\cO _X(-p)\otimes _{\cO _X}\cB _n]=[\cO _p\otimes
_{\cO _X}\cB _n]$$ in $K_0(X)$.  By \cite [Lemma
4.1]{Ng}, the length of $\cO _p\otimes _{\cO _X}\cB _n\in \mod X$
is $n+1$, so $[\cO _p\otimes _{\cO _X}\cB _n]=\sum _{i=0}^n[\cO
_{p_i}]$ in $K_0(X)$ for some closed points $p_0, \cdots , p_n\in
X$.  In particular,
\begin{align*}
(f^*\cO _X , (f^*\cO _p)(n)) & =(\cO _X ,
f_*((f^*\cO _p)(n)))-(\cO _X, R^1f_*((f^*\cO _p)(n))) \\
& =\sum _{i=0}^n(\cO _X , \cO _{p_i})
\\
& =n+1
\end{align*}
by Lemmas \ref{lem.prf} and \ref{lem.curve}.
\end{pf}

The next result says that ``fibers do not meet'' and ``a fiber and a
section meet exactly once''.

\begin{theorem}
\label{thm.int.thy}
Let $\PP (\cE)$ be a quantum ruled surface over $X$, and $p$ and $q$
closed points of $X$.
Define
$$
H:=[\cO_{\PP(\cE)}]-[\cO_{\PP(\cE)}(-1)] \in K_0(\PP(\cE)).
$$
Then
\begin{align*}
& f^*\cO _p\cdot f^*\cO _q=0. \\
& f^*\cO _p\cdot H=1. \\
& H\cdot f^*\cO _p=1 \\
& H\cdot H=\deg (pr_{2*}\cE ).
\end{align*}
\end{theorem}
\begin{pf}
(1) By Lemmas \ref{lem.prf} (2), \ref{lem.curve}, and \ref {lem.ed},
$$(f^*\cO _p, f^*\cO _q)
 =(\cO _p, f_*(f^*\cO _q))-(\cO _p, R^1f_*(f^*\cO _q))
 =(\cO _p, \cO _q)
 =0.
$$

(2)
By Lemmas \ref{lem.prf} (3) and \ref{lem.curve}
\begin{align*}
(f^*\cO _p, H) & =(f^*\cO _p,
f^*\cO _X) - (f^*\cO _p, (f^*\cO _X)(-1))
\\
& =(\cO _p, \cO _X\otimes _{\cO _X}\cB _0) - (\cO _p,
\cO _X\otimes _{\cO _X}\cB _{-1})
\\
& =(\cO _p, \cO _X)
\\
&=-1.
\end{align*}

(3)  By Lemma \ref {lem.ed},
\begin{align*}
(H, f^*\cO _p) & = (f^*\cO _X,
f^*\cO _p) - ((f^*\cO _X)(-1), f^*\cO _p)
\\
& = (f^*\cO _X, f^*\cO _p) - (f^*\cO _X, (f^*\cO _p)(1))
\\
& =1-2
\\
& =-1.
\end{align*}

(4) Let $g$ be the genus of $X$.  Since
$pr _{2*}\cE$ is a locally free module of rank 2 on a smooth curve,
there is an exact sequence
$$0\to \cL _1\to pr _{2*}\cE \to \cL _2\to 0$$ in
which $\cL _1$ and $\cL _2$ are line bundles
\cite [Chapter V, Exercise 2.3 (a)]{H}.
Thus $[pr_{2*}\cE]=[\cL_1]+[\cL_2]$
and $\deg (pr_{2*}\cE)=\deg \cL _1+\deg \cL _2$ \cite
[Chapter II, Exercise
6.12]{H}. By Riemann-Roch, and Lemmas \ref{lem.prf} (3) and
\ref{lem.curve}, we have
\begin{align*}
(H, H)
= & (f^*\cO _X, f^*\cO _X)
-((f^*\cO _X)(-1), f^*\cO _X)  \\
& \; \; \; \; \; \; \; \; \; \; - (f^*\cO _X,
(f^*\cO _X)(-1)) + ((f^*\cO _X)(-1), (f^*\cO
_X)(-1))
\\
= & (\cO _X, \cO _X\otimes _{\cO _X}\cB _0) -(\cO _X,
\cO _X\otimes _{\cO _X}\cB _1)
- (\cO _X, \cO _X\otimes_{\cO _X}\cB _{-1})
\\
& \qquad \qquad
+ (\cO _X, \cO _X\otimes _{\cO _X}\cB _0)
\\
= & (1-g) - (1-g+\deg \cL _1)-(1-g+\deg \cL _2) +(1-g)\\
= & -\deg \cL _1-\deg \cL _2
\\
= & -\deg (pr_{2*}\cE ).
\end{align*}
\end{pf}

The previous result is the same as for the commutative
case \cite[Propositions V.2.8 and V.2.9]{H}.
When $\cE$ is ``normalized", it is
reasonable to call $H$ a K-theoretic section; in this case
the integer $e:=-H\cdot H =-\deg(pr_{2*}\cE)$ should play an
important role in the classification of quantum ruled
surfaces analogous to the role that $e$ plays in the commutative
theory.

The rank of a coherent $\cO_X$-module gives a group
homomorphism $\rank:K_0(X)
\to \ZZ$. Recall that $K_0(\PP(\cE))=K_0(X) \oplus K_0(X)T$,
where the $K_0(X)$
on the right hand side is really $f^*K_0(X)$. We define $\rank:K_0(\PP(\cE))
\to \ZZ$ by $\rank(a+bT)=\rank(a)+\rank(b)$ where $a,b \in K_0(X)$. Then
$\rank(\cO_{\PP(\cE)}(n))=1$ for all $n$, so $\rank(H)=0$, and also $\rank(f^*\cO_p
(n)) =0$ for all $n$ and all $p \in X$. We define $F^1K_0(\PP(\cE))$ to be the
kernel of the rank function and think of it as the subgroup of $K_0(\PP(\cE))$
generated by those $\PP(\cE)$-modules whose ``support'' is
of codimension $\ge 1$.

The kernel of $\rank:K_0(X) \to \ZZ$ is $F^1K_0(X)$, the subgroup
generated by the $\cO_X$-modules having support of codimension $\ge
1$, and this subgroup is isomorphic to the Picard group $\Pic X$.
We have $F^1K_0(\PP(\cE)) = F^1K_0(X) \oplus F^1K_0(X).H \oplus \ZZ H$.
Using the results already obtained,
one sees that $F^1K_0(X).H$ is contained in both the left and right
radicals of the Euler form restricted to $F^1K_0(\PP(\cE))$, so the Euler
form induces a $\ZZ$-valued bilinear form on the quotient
$F^1K_0(\PP(\cE))/F^1K_0(X).H$. There are good reasons to think of this
quotient as playing the role of the Picard group for $\PP(\cE)$, and if one does
this one has $\Pic \PP(\cE) \cong f^* \Pic X \oplus \ZZ.H$ just as in
the commutative case (cf. \cite[Ch. V, Prop.2.3]{H}); here $f^*\Pic X$ denotes the
image of $F^1K_0(X)$ under the map $f^*:K_0(X) \to K_0(\PP(\cE))$.

Modding out all the radical for the Euler form on $F^1K_0(\PP(\cE))$ produces
an analogue of the Neron-Severi group, and this is a free abelian group
of rank two with basis
$[f^*\cO_p]$ and $H$, and the pairing on this is given by Theorem \ref{thm.int.thy}.
Again, one has a good analogue of the commutative theory.

\ifx\undefined\bysame
\newcommand{\bysame}{\leavevmode\hbox to3em{\hrulefill}\,}
\fi

\end{document}